\documentclass{article}
\usepackage{graphicx} 
\usepackage{amsmath}
\usepackage{amsfonts}
\usepackage{amsthm}
\usepackage{url}
\usepackage{color}
\newtheorem{theorem}{Theorem}[section]

\newtheorem{remark}[theorem]{Remark}

\usepackage{graphicx}
\usepackage{hyperref}
\usepackage{subfig}
\usepackage{color}
\usepackage{diagbox}
\usepackage{multirow}
\usepackage{booktabs}
\usepackage{threeparttable}        
\usepackage[numbers,sort&compress]{natbib}
\usepackage{CJKutf8} 
\usepackage{makecell}
\usepackage{cleveref}
\usepackage[lined,boxed,linesnumbered,ruled]{algorithm2e}
\usepackage{bm}

\begin{document}
\date{}
\title{High-Precision Phase-Shift Transferable Neural Networks for High-Frequency Function Approximation and PDE Solution}
\author{~~Xuyang Gao
\thanks{School of Mathematics and Science, Harbin Normal University, Harbin, 150025, China.},
~~Liang Chen\thanks{Department of Mathematics, Jiujiang University, Jiujiang, 332000, P. R. China.},
~~Minqiang Xu\thanks{Corresponding Author. School of Mathematical Sciences, Zhejiang University of Technology, Hangzhou, 310023, China(mqxu@zjut.edu.cn).},
~~Jing Niu\thanks{School of Mathematics and Science, Harbin Normal University, Harbin, 150025, China.}
}
\maketitle

\textbf{Abstract:} 
Neural network-based methods have emerged as a promising paradigm for scientific computing, yet they face critical bottlenecks in high-frequency function approximation and partial differential equation (PDE) solving. A fundamental limitation is the inherent spectral bias of neural networks, which prioritizes low-frequency components and hinders the accurate capture of sharp, oscillatory details. 

To address these challenges, we propose the Phase-Shift Transferable Neural Network (PhaseTNN) framework. By integrating phase-shift frequency decomposition with TransNet’s efficient least-squares training, PhaseTNN eliminates the need for manual shape parameter tuning and decomposes intractable high-frequency problems into solvable low-frequency subproblems. We present two complementary architectures: Parallel PhaseTNN (PPTNN), designed for parallel frequency learning, and Coupled PhaseTNN (CPTNN), tailored for simultaneous multi-frequency representation.

Extensive numerical experiments demonstrate that PhaseTNN achieves near-machine-precision accuracy for high-frequency approximation.  For high-frequency PDEs (variable-coefficient elliptic equation, linear/nonlinear Helmholtz, Wave equations, diffusion interface problem), CPTNN consistently outperforms baselines by up to three orders of magnitude in accuracy while delivering orders-of-magnitude faster computation.
These results confirm PhaseTNN’s superior robustness, accuracy, and efficiency for high-frequency scientific computing tasks.



\section{Introduction}
In recent years, neural networks have achieved remarkable progress across diverse fields, including computer vision, speech recognition, natural language processing, and scientific computing \cite{lecunDeepLearning2015,ronnebergerUnetConvolutionalNetworks2015,heDeepResidualLearning2016,mildenhallNerfRepresentingScenes2021,collobertUnifiedArchitectureNatural2008,chowdharyNaturalLanguageProcessing2020,brownStatisticalApproachMachine1990a,vaswaniAttentionAllYou2017,devlinBertPretrainingDeep2019,brownLanguageModelsAre2020,hanSolvingHighdimensionalPartial2018,hanDeepLearningbasedNumerical2017,raissiPhysicsinformedNeuralNetworks2019,sirignanoDGMDeepLearning2018,yuDeepRitzMethod2018,zangWeakAdversarialNetworks2020}. Against this backdrop, neural network-based methods for solving partial differential equations (PDEs) have emerged as a transformative paradigm, garnering widespread attention in the scientific computing community. Representative approaches include Physics-Informed Neural Networks (PINNs)
\cite{raissiPhysicsinformedNeuralNetworks2019}, the Deep Galerkin Method (DGM) \cite{sirignanoDGMDeepLearning2018}, the Deep Ritz Method (DRM) \cite{yuDeepRitzMethod2018}, and Weak Adversarial Networks (WAN) \cite{zangWeakAdversarialNetworks2020}.

However, these methods face challenges that impede their practical utility in scientific computing. On one hand, neural network training is prone to stagnation in local minima \cite{baldiNeuralNetworksPrincipal1989,krishnapriyanCharacterizingPossibleFailure2021}, leading to premature convergence of the loss function. This not only degrades the final accuracy of the model but also introduces failure modes inherent to deep neural networks (DNNs) \cite{krishnapriyanCharacterizingPossibleFailure2021}. On the other hand, the innate spectral bias of neural networks undermines their ability to capture sharp variations in PDE solutions, a phenomenon well-documented in the literature \cite{johnxuFrequencyPrincipleFourier2020,rahamanSpectralBiasNeural2019,luoUpperLimitDecaying2022,xuOverviewFrequencyPrinciple2025,xuUnderstandingOvercomingSpectral2025a}. This spectral bias creates a critical bottleneck for learning high-frequency components in PDEs, rendering high-frequency function approximation and PDE solving among the most demanding scenarios for neural network-based methods. Consequently, achieving high-precision solutions remains a fundamental prerequisite and core challenge for deploying such methods in high-performance scientific computing.

To address the aforementioned challenges, existing research has primarily advanced along two distinct paradigms. The first paradigm centers on enhancing gradient-descent-based optimization frameworks, with specific strategies including advances in optimization algorithms  
\cite{chiuCANPINNFastPhysicsinformed2022,mullerAchievingHighAccuracy2023,bottouOptimizationMethodsLargescale2018,wangGradientAlignmentPhysicsinformed2025,urbanUnveilingOptimizationProcess2025,kiyaniOptimizingOptimizerPhysicsinformed2025}, architectural innovations \cite{jagtapExtendedPhysicsinformedNeural2020,moseleyFiniteBasisPhysicsinformed2023,ainsworthGalerkinNeuralNetwork2022,liuKanKolmogorovarnoldNetworks2024,shuklaComprehensiveFAIRComparison2024}, loss function term weighting schemes \cite{wangUnderstandingMitigatingGradient2021,wangWhenWhyPINNs2022}, and improvements to training protocols \cite{howardStackedNetworksImprove2023,kopanicakovaEnhancingTrainingPhysicsinformed2024}. Despite these developments, achieving high-precision predictions remains a persistent, severe challenge. A recently proposed multi-stage neural network (MSNN) approach
\cite{wangMultistageNeuralNetworks2024} employs a hierarchical training strategy that iteratively refines residuals, enabling predictions to approach machine precision. While this method has achieved notable progress in approximation accuracy, it still faces fundamental limitations in training efficiency and generalization performance. Consequently, gradient-based methods currently encounter inherent difficulties in balancing solution accuracy and computational cost.

The second paradigm seeks to fundamentally restructure the training pipeline to overcome the optimization challenges and high computational overhead inherent to deep neural networks (DNNs) \cite{hornikMultilayerFeedforwardNetworks1989}. Representative methods in this category include the Extreme Learning Machine (ELM) \cite{huangExtremeLearningMachine2006} and the Random Feature Method (RFM) \cite{chenBridgingTraditionalMachine2022}. Their core design is to construct a shallow neural network: hidden-layer parameters are randomly initialized and fixed, leaving only the output layer trainable. This formulation transforms the learning problem into an efficiently solvable linear least-squares problem. Inspired by physics-informed neural networks (PINNs) and ELM, Dwivedi et al. \cite{dwivediPhysicsInformedExtreme2020} proposed the Physics-Informed Extreme Learning Machine (PIELM), which combines the unsupervised learning characteristics of PINNs with the high computational efficiency of ELM. The local extreme learning machine (locELM) \cite{dongLocalExtremeLearning2021} further integrates ELM with domain decomposition techniques, significantly enhancing its capability to solve complex problems. Overall, compared with traditional gradient-based neural networks, these approaches substantially reduce computational costs and often achieve high accuracy for PDE solving \cite{chenBridgingTraditionalMachine2022,dongLocalExtremeLearning2021,dongComputingHyperparameterExtreme2022,chenRandomFeatureMethod2023,yangNovelImprovedExtreme2018,niNumericalComputationPartial2023}. Nevertheless, ELM and RFM lack rigorous, theoretically justified guidelines for selecting hidden-layer parameters; parameter tuning relies heavily on empirical experience or manual adjustment, limiting their reliability and scalability in large-scale, high-frequency problems.

The recently proposed \textbf{Transferable Neural Network (TransNet)} \cite{zhangTransferableNeuralNetworks2024}  effectively addresses this scalability limitation. TransNet constructs a single-hidden-layer network where the weights and biases of hidden neurons are reparameterized via location and shape parameters. By tuning shape parameters via auxiliary functions, it generates uniformly distributed, transferable neurons. Numerical results confirm that shape parameter selection works effectively for low-frequency problems, enabling TransNet to achieve high computational accuracy and efficiency.

While TransNet performs well for low-frequency problems, its extension to high-frequency regimes remains a critical challenge. A key limitation lies in its heavy reliance on the \textbf{shape parameter}, which governs the steepness of hidden neurons and dictates the expressivity of the approximation space. Existing tuning strategies work effectively for low-frequency scenarios but fail in high-frequency settings, suffering from three core drawbacks:
\begin{itemize}
    \item  Gaussian Random Field (GRF)-based tuning methods \cite{zhangTransferableNeuralNetworks2024} introduce an additional correlation-length hyperparameter, whose optimal value is tightly coupled to the behavior of the underlying PDE solution. Improper selection of this length scale compromises the GRF’s ability to capture the true solution structure, leading to misaligned shape parameters. 
    \item  The optimal shape parameter depends strongly on network scale, even for the same problem, rendering the tuning process computationally expensive. Although an empirical extrapolation strategy \cite{luMultipleTransferableNeural2025} mitigates this cost, its accuracy is constrained by the generality of the empirical constant. 
   \item A fixed shape parameter cannot simultaneously capture low-frequency backgrounds and high-frequency details required for high-precision solutions.
\end{itemize}

These observations reveal a fundamental limitation of TransNet: its insufficient approximation capability for high-frequency problems, rooted in \textbf{a parameterization strategy that fails to align with the multi-scale frequency characteristics of PDE solutions}. To overcome this challenge, we require a framework that can effectively capture the features of rapidly oscillatory solutions while preserving computational efficiency.

Among existing efforts to address high-frequency challenges, the phase-shift deep neural network \cite{caiPhaseShiftDeep2020} provides particularly critical inspiration: it decomposes challenging high-frequency problems into multiple tractable low-frequency subproblems via frequency-domain filtering and phase shifting, thereby significantly improving the DNN’s ability to fit high-frequency components. However, PhaseDNN offers a profound conceptual foundation but relies on gradient descent training, which introduces significant optimization challenges and high computational costs. In contrast, TransNet leverages efficient least-squares solvers to deliver rapid training. We therefore aim to combine the strengths of both approaches: a method tailored for high-frequency solutions while preserving computational efficiency.

To this end, we propose the \textbf{Phase Shift Transferable Neural Network (PhaseTNN)} framework, encompassing two architectures: the Parallel PhaseTNN (PPTNN) and the Coupled PhaseTNN (CPTNN). PPTNN decomposes the target problem into distinct frequency bands via frequency-domain filtering, then transforms each component into a tractable low-frequency subproblem via phase shifting, enabling these components to be learned in parallel. Instead of explicitly decomposing training data, CPTNN constructs phase-modulated neural basis functions to simultaneously represent multiple frequency components within a single model, enabling efficient learning for high-dimensional and large-scale problems.

The main contributions of our work are summarized as follows:
\begin{itemize}
    \item \textbf{An improved filtering method}: Compared with the filtering method used in PhaseDNN \cite{caiPhaseShiftDeep2020}, our method achieves significantly improved filtering accuracy, approaching double-precision machine precision limits.   
    \item \textbf{Enhanced robustness and usability}: Compared with TransNet methods in \cite{zhangTransferableNeuralNetworks2024, luMultipleTransferableNeural2025}, our method avoids shape parameter tuning and possesses stronger adaptability. For high-frequency problems, conventional TransNets are highly sensitive to the shape parameter, rendering the strategies in \cite{zhangTransferableNeuralNetworks2024, luMultipleTransferableNeural2025} ineffective; whereas our method achieves robust and superior performance across various high-frequency functions.    
    \item \textbf{High accuracy and exceptional computational efficiency}: Numerical experiments confirm that our proposed method significantly outperforms many existing techniques in computational speed, while achieving high-precision approximation of high-frequency functions and accurate solution of PDEs with high-frequency components.
\end{itemize}

The remainder of the paper is organized as follows. In \cref{sec:related work}, we establish the transferable neural feature space. \cref{sec:methodology approximation} presents the PhaseTNN framework, detailing two distinct architectures, the PPTNN and the CPTNN, and their application to function approximation. Subsequently, \cref{sec:methodology solving} introduces the application of these architectures to solving differential equations. Comprehensive numerical experiments demonstrating the effectiveness of our method are presented in \cref{sec:numerical results}. Finally, \cref{sec:conclusion} concludes the paper with a summary of key findings.

\section{Construction of a Transferable Neural Feature Space for low-frequency functions}
\label{sec:related work}
Inspired by \cite{zhangTransferableNeuralNetworks2024}, we  consider a single-hidden-layer fully-connected neural network, denoted by
\begin{equation}
T_{\text{NN}}(\boldsymbol{x},\boldsymbol{\alpha}) := \sum_{m=1}^M \alpha_m \sigma(\boldsymbol{w}_m^\top \boldsymbol{x } + b_m) + \alpha_0,
\label{eq:single_layer_nn}
\end{equation}
where \( \boldsymbol{x} \in \Omega \subset \mathbb{R}^d \), $M$ is the number of hidden neurons,  $\boldsymbol{w}_m = (w_{m,1}, \dots, w_{m,d})^\top$ and $b_m$ are 
are the weight vector and bias of the $m$-th hidden neuron, the vector $\boldsymbol{\alpha} = (\alpha_0, \alpha_1, \dots, \alpha_M)^\top$ collects the output-layer weights and bias, and $\sigma(\cdot)$ denotes the activation function.

In \cite{zhangTransferableNeuralNetworks2024}, the authors establish a transferable neural feature space $\mathcal{P}_{\text{NN}}$ defined as
\begin{equation}
\mathcal{P}_{\text{NN}} = \operatorname{span}\left\{1,\ \sigma(\boldsymbol{w}_1^\top \boldsymbol{x} + b_1),\ \dots,\ \sigma(\boldsymbol{w}_M^\top \boldsymbol{x} + b_M)\right\},
\label{eq:neural_feature_space}
\end{equation}
To obtain uniformly distributed basis functions over  $\mathcal{P}_{\text{NN}}$, the key strategy adopted in \cite{zhangTransferableNeuralNetworks2024} is to reparameterize the hidden-layer parameters $(\boldsymbol{w}_m, b_m)$ in terms of location parameters $(\boldsymbol{a}_m, r_m)$ and a shape parameter $\gamma_m$, such that
\[\sigma(\boldsymbol{w}_m^\top \boldsymbol{x}+b_m)=\sigma(\gamma_m(\boldsymbol{a}_m^\top \boldsymbol{x}+r_m)),\]
with $\|\boldsymbol{a}_m\|_2 = 1$ and $\gamma_m \geq 0$. Their equivalence  between the two parameterizations is given by
\[
\begin{cases}
\boldsymbol{w}_{m} = \gamma_{m} \boldsymbol{a}_{m} \\
 b_m = \gamma_m r_m
\end{cases}
\iff
\begin{cases}
\boldsymbol{a}_m = \frac{\boldsymbol{w}_m}{\|\boldsymbol{w}_m\|_2} \\
r_m = \frac{b_m}{\|\boldsymbol{w}_m\|_2} \\
\gamma_m = \|\boldsymbol{w}_m\|_2
\end{cases}.
\]
The complete sampling procedure for the location parameters is expressed as follows:
\begin{equation}
\boldsymbol{a}_m = \frac{X_m}{\| X_m \|_2} \quad \text{and} \quad r_m = U_m, \ m = 1, \dots, M.
\label{eq:L_neural_feature_space_am}
\end{equation}
where $X_m$ are independent and identically distributed (i.i.d.) standard Gaussian random vectors, $U_m$ follows i.i.d. uniform random variables on $[0, 1]$ and the optimal shape parameters  $\gamma_m$ is obtained via pre-training. 

While the approximation space  $\mathcal{P}_{\text{NN}}$ achieves satisfactory performance for low-frequency functions, the pre-trained optimal parameters fail to generalize to high-frequency functions. Moreover, exhaustive enumeration for optimal shape parameter selection during training incurs prohibitive computational costs. To address these limitations, inspired by \cite{caiPhaseShiftDeep2020}, we propose to extract high-frequency components from target functions via filtering, and transform these components into tractable low-frequency components via a translation (phase-shift) transformation. In this framework, we only require   $\mathcal{P}_{\text{NN}}$ to learn the resulting low-frequency components, eliminating the need for costly shape parameter tuning. Consequently, the shape parameters $\gamma_m$ can be chosen randomly, and we construct the low-frequency transferable feature space
$\mathcal{P}_{\text{NN}}^{L}$ as
\begin{equation}
\mathcal{P}_{\text{NN}}^L = \operatorname{span}\left\{1,\ \sigma\left(\gamma_1(\boldsymbol{a}_1^\top \boldsymbol{x} + r_1)\right), \dots,\ \sigma\left(\gamma_M(\boldsymbol{a}_M^\top \boldsymbol{x} + r_M)\right)\right\},
\label{eq:L_neural_feature_space}
\end{equation}
where $\boldsymbol{a}_m$ and $r_m$ are obtained via \eqref{eq:L_neural_feature_space_am}, and $\gamma_m$ is fixed to a small constant—specifically, $\gamma_m=2$ throughout this work.

\begin{remark}
The filtering operation effectively converts high-frequency problems into equivalent low-frequency problems for solution. Therefore, in our method, \( \gamma_m \)  is simply fixed to a shape parameter that is nearly optimal for most low-frequency problems, thereby avoiding the laborious shape parameter tuning process required in \cite{zhangTransferableNeuralNetworks2024}.
\end{remark}

\section{Phase-Shift Transferable Neural Networks for High-Frequency Approximation}
\label{sec:methodology approximation} 

\subsection{An improved filtering method}
In this subsection, we address two core challenges: (\textbf{i}) how to apply frequency-domain filtering to a function  $f$ that is only defined and observable on \([-1, 1]\), (\textbf{ii}) how to integrate this filtering operation into numerical approximations over the entire interval.

We first define the Fourier transform and its inverse for a function \( f(x) \) by
\begin{equation}
\mathcal{F}[f](k) = \int_{-\infty}^{+\infty} f(x) e^{-2\pi \mathrm{i}kx} \, dx, \quad \mathcal{F}^{-1}[\hat{f}](x) =  \int_{-\infty}^{+\infty} \hat{f}(k) e^{2\pi \mathrm{i}kx} \, dk. 
\end{equation}
For a given frequency increment $\Delta k$, we assume that for some integer $K > 0$, the support of $\operatorname{supp} \hat{f}(k)$ is contained in
$ [-K\Delta k, K\Delta k]$, i.e.,
\[
\operatorname{supp} \hat{f}(k) \subset [-K\Delta k, K\Delta k].
\]
For any positive integer $N$, we denote $\mathcal{I}_{N} = \{-N, \dots, N\}$ and $\mathcal{I}_{N}^+ = \{1,2, \dots, N\}$. 
We first construct a mesh for the interval $[-K\Delta k, K\Delta k]$ by
\begin{equation}
\kappa_j = j\Delta k,~j \in \mathcal{I}_{K}
\end{equation}
Let $\phi(k)$ denote the characteristic function of the interval
$\left[-\frac{1}{2}, \frac{1}{2}\right]$, i.e.,
$\phi(k) = \chi_{\left[-\frac{1}{2}, \frac{1}{2}\right]}(k)$, with inverse Fourier transform $\phi^\vee(x) =  \frac{\sin \pi x}{\pi x}$ (the sinc function). We then define the scaled basis functions
$$\phi_j(k) = \phi\!\left(\frac{k - \kappa_j}{\Delta k}\right),~j \in \mathcal{I}_{K},$$
which satisfy the partition-of-unity identity
\begin{equation}
1 = \sum_{j=-K}^K \phi_j(k), \ k \in [-K\Delta k, K\Delta k].
\end{equation}
Define the frequency-shifted basis function via inverse Fourier transform:
\begin{equation}
\phi_j^\vee(x):=\mathcal{F}^{-1}[\phi_j](x)=\Delta ke^{2\pi \mathrm{i}\kappa_jx} \phi^\vee(\Delta k x), x \in  \mathbb{R}.
\end{equation}
Accordingly, the target function  $f(x)$ decomposes in the frequency domain as
\begin{equation}
\hat{f}(k) = \sum_{j=-K}^K \phi_j(k) \hat{f}(k) := \sum_{j=-K}^K \hat{f}_j(k). 
\end{equation}
The corresponding spatial-domain decomposition of $f$ reads
\begin{equation}
f(x) = \sum_{j=-K}^K f_j(x),~\text{where}~ f_j(x) = \mathcal{F}^{-1}[\hat{f}_j](x),j \in \mathcal{I}_{K}.
\end{equation}
For any $x_i\in[-1,1]$, the pointwise value $f_j(x_i)$ is computed as
\begin{equation}
\begin{split}
f_j(x_i)& = \mathcal{F}^{-1}[\hat{f}_j](x_i)=\mathcal{F}^{-1}[\phi_j(k) \hat{f}(k)](x_i)=\mathcal{F}^{-1}[\phi_j(k)]* \mathcal{F}^{-1}[\hat{f}(k)](x_i) \\
&=(\phi_j^\vee * f)(x_i) = \int_{-\infty}^{+\infty} \phi_j^\vee(x_i - s) f(s) ds,j \in \mathcal{I}_{K}.
\end{split}
\label{jfjuanji}
\end{equation}

If $f$ is only known on \([-1, 1]\),  how do we evaluate the integral in \eqref{jfjuanji}? Our goal is to construct an auxiliary function  $F$ that coincides exactly with $f$ over the interval \([-1, 1]\), exhibits at least
$(Q-1)$-order continuous differentiability at the endpoints \(x=-1,1\), and decays rapidly to zero outside this interval. To this end, we employ a higher-order exponential decay factor. Specifically, the function $F$ is constructed via two distinct strategies. 

Case 1: Explicit expression of $f$ on \([-1, 1]\) is available,
\begin{equation}
F(x) = 
\begin{cases}
f(x)\text{e}^{-10(x+1)^Q}, & x \in (-\infty, -1),\\
~~~~f(x), & x \in [-1, 1],\\
f(x)\text{e}^{-10(x-1)^Q},  & x \in (1, +\infty).\\
\end{cases}
\label{eq:hf1}
\end{equation}

Case 2: Only finitely many sampling points of $f$ on \([-1, 1]\) are known. We construct polynomial extensions on both sides of $f$ on \([-1, 1]\) to define $F$:
\begin{equation}
F(x) = 
\begin{cases}
\left(\sum\limits_{q=0}^{Q-1} \frac{f^{(q)}(-1)}{q!}(x+1)^q\right)\text{e}^{-10(x+1)^Q}, & x \in (-\infty, -1),\\
~~~~~~~f(x), & x \in [-1, 1],\\
~\left(\sum\limits_{q=0}^{Q-1} \frac{f^{(q)}(1)}{q!}(x-1)^q\right)\text{e}^{-10(x-1)^Q},  & x \in (1, +\infty).\\
\end{cases}
\label{eq:hf2}
\end{equation}
where the values of $f^{(q)}(-1)$ and $f^{(q)}(1)$ for $q=0,1,\ldots,Q-1$
are approximated from samples of $f$ on \([-1, 1]\).

Once the data are prepared, we can approximate  \(f_j(x)\) via the trapezoidal rule over the domain $[-C, C]$:
\begin{equation}
f_j(x)\approx\int_{-C}^{C} \phi_j^\vee(x - s) F(s) ds\approx \sum_{s=1}^{N_s} \phi_j^\vee(x - x_s) F(x_s) \Delta \tilde{x}_s :=\tilde{f}_j(x),
\label{eq:traningdate}
\end{equation}
where $\{x_s\}_{s=1}^{N_s}$ denote numerical integration nodes with grid spacing $\Delta x = \frac{2C}{N_s - 1}$, and $C$ is chosen such that the function $F(x_s)$ is sufficiently small outside  $(-C, C)$, and $\Delta \tilde{x}_s$ are given by:
\begin{equation}
\Delta \tilde{x}_s = 
\begin{cases} 
\frac{\Delta x}{2}, & \text{if } s = 1 \text{ or } s = N_s, \\ 
\Delta x, & \text{if } 1 < s < N_s. 
\end{cases}
\end{equation}

Consequently, the overall approximation $\tilde{f}$ of $f$ reconstructed as
\begin{equation}
\tilde{f}=\sum_{j=-K}^{K}\tilde{f}_j(x),\quad x\in[-1,1].
\label{eq:appromatingf}
\end{equation}

It is important to note that our method differs from \cite{caiPhaseShiftDeep2020}: while \cite{caiPhaseShiftDeep2020} only considers scenarios where \(f(x)\) has an explicit expression, our work addresses scenarios where $f$ is known over \([-1, 1]\) (or only finite sampled data is available). Furthermore,  when using the filtering method in \cite{caiPhaseShiftDeep2020} to approximate \(f(x)\), our method significantly outperforms  \cite{caiPhaseShiftDeep2020}  in terms of approximation accuracy and can even achieve machine precision. Detailed comparisons are provided in Example \ref{sec:filtering accuracy} of Section \ref{sec:numerical results}.

\begin{remark}
Compared with the method in  \cite{caiPhaseShiftDeep2020}, our approach offers a distinct advantage through the introduction of a decay factor. This factor guarantees rapid decay of the integrand, enabling high-accuracy approximation without requiring an excessively large integration domain or an inordinate number of sampling points. Throughout this work, the parameter \( C \) is set to 5.
\end{remark}
\begin{remark}
Suppose only finitely many values of $f$ are available on the interval $[-1,1]$. We employ the PPR method \cite{zhangNewFiniteElement2005,xuNovelClassHessian2025}  to recover the derivatives of  $f$  at the endpoints $x=-1$ and $x=1$. For each \(z\in\{-1, 1\}\), we first construct a local stencil $\mathcal{L}_z$ around  $z$ that contains at least $Q$ points, and then determine the polynomial \(p_{z} \in \mathbb{P}_{Q-1}(\mathcal{L}_z)\) that solves the following least-squares minimization problem:
\[
p_z = \arg\min_{p \in \mathbb{P}_{Q-1}(\mathcal{L}_z)} \sum_{\tilde{z} \in \mathcal{L}_z } \left| p(\tilde{z}) - v_h(\tilde{z}) \right|^2.
\]
The $q$-th derivative of $f$ at $z$ is then approximated by evaluating the $q$-th derivative of the polynomial  $p_{z}$ at $z$:
\[
f^{(q)}(z) = \frac{\partial^q p_z}{\partial x^q} \bigg|_{x=z}, \quad q = 0, 1, \dots, Q-1.
\]
\end{remark}

\subsection{A Parallel Phase-Shift Transferable Neural Network (PPTNN) for the 1D Case}
In this subsection, we develop a framework for high-accuracy approximation of high-frequency functions, leveraging the transferable neural feature space established in Section \ref{sec:related work} to approximate each component \( \tilde{f}_j(x) \) in Eq. \eqref{eq:appromatingf}. Our primary contribution is the proposal of a \textbf{Parallel Phase-Shift Transferable Neural Network (PPTNN)}, which leverages phase-shift frequency decomposition to enable efficient learning of each frequency component.

Following the strategy in \cite{caiPhaseShiftDeep2020}, we implement the learning process for each individual frequency component \( \tilde{f}_j(x) \), whose frequency spectrum is confined to the interval \( [\kappa_j - \frac{\Delta k}{2}, \kappa_j + \frac{\Delta k}{2}] \).  It is a well-established result in harmonic analysis that low-frequency components of a loss function converge significantly faster than their high-frequency counterparts. Exploiting this property, a simple phase shift translates the spectrum of \( \tilde{f}_j(x) \) into the baseband interval \( [-\frac{\Delta k}{2}, \frac{\Delta k}{2}] \), thereby enabling efficient training of a Transferable Neural Network \( T_j(x) \) with only a few training epochs.

Specifically, given that the frequency support of
 \( \hat{f}_j(k) \) is \( [\kappa_j - \frac{\Delta k}{2}, \kappa_j + \frac{\Delta k}{2}] \), the shifted spectrum \( \hat{f}_j(k + \kappa_j) \) is supported in \( [-\frac{\Delta k}{2}, \frac{\Delta k}{2}] \). Its inverse Fourier transform is defined as 
\begin{equation}
f_j^{\text{shift}}(x) = \mathcal{F}^{-1} \left[ \hat{f}_j(k + \kappa_j) \right] (x), j \in \mathcal{I}_{K}.
\label{eq:2.12}
\end{equation}
Subsequently, we focus on approximating
$f_j^{\text{shift}}(x)$ in $\mathcal{P}_{\text{NN}}^L$ using a transferable neural network $ T_j(x,\boldsymbol{\alpha}_j)$, defined as
$$ T_j(x,\boldsymbol{\alpha}_j)=\sum_{m=0}^N\alpha_{j,m} \psi_{j,m}(x),j \in \mathcal{I}_{K}.$$
To construct the training dataset, we utilize the fundamental frequency-shift identity:
\begin{equation}
f_j^{\text{shift}}(x) = e^{-2\pi i\kappa_j x} f_j(x), j \in \mathcal{I}_{K}.
\label{eq:frequency shift}
\end{equation}
Recalling the approximation  $\tilde{f}_j(x)$ in \eqref{eq:traningdate}, the training data $\{x_{\xi},f_j^{\text{shift}}(x_{\xi})\}_{\xi=1}^{N_{\xi}}$ are obtained via numerical quadrature:
\begin{align}
f_j^{\text{shift}}(x_\xi) \approx e^{-2\pi i \kappa_j x_\xi} \sum_{s=1}^{N_s} \phi_j^\vee(x_{\xi} - x_s) F(x_s) \Delta \tilde{x}_s :=\tilde{f}_j^{\text{shift}}(x_{\xi}), \quad \xi \in \mathcal{I}_{N_{\xi}}^+,
\label{eq:shfittraningdate}
\end{align}
where
\begin{equation}
\label{eq:trapz_weights}
\Delta \tilde{x}_s = 
\begin{cases} 
\frac{\Delta x}{2}, & \text{if } s = 1 \text{ or } s = N_s, \\ 
\Delta x, & \text{if } 1 < s < N_s. 
\end{cases}
\end{equation}
The network $ T_j(x,\boldsymbol{\alpha}_j)$ is trained by minimizing the mean-squared error (MSE) loss function:
\begin{equation}
L_j(\boldsymbol{\alpha}_j) = \sum_{\xi=1}^{N_{\xi}} \left( f_j^{\text{shift}}(x_{\xi}) - T_j(x_{\xi}, \boldsymbol{\alpha}_j) \right)^2,j \in \mathcal{I}_{K}.
\label{determT}
\end{equation}
Combining \eqref{eq:frequency shift} and the definition of \eqref{determT}, we obtain the approximation for $\tilde{f}_j(x)$:
\begin{equation}
 \tilde{f}_j(x)= e^{2\pi i\kappa_j x}T_j(x, \boldsymbol{\alpha}_j). 
\label{eq:approfj}
\end{equation}
Finally, the overall approximation \(f(x)\) is reconstructed by summing the component-wise approximations: 
\begin{equation}
f(x)\approx \sum_{j=-K}^{K}\tilde{f}_j(x):=\tilde{f}(x), \quad\text{for}\quad x\in(-1,1).
\end{equation}

\begin{remark}
In practice, the proposed PPTNN framework does not require training of all $2K+1$ frequency components \( f_j(x) \).  Instead, we prioritize learning dominant frequency components that carry the majority of signal energy, while extremely weak components can be safely approximated as zero. This adaptive learning strategy directs the network to focus on problem-relevant frequency content, thereby significantly reducing the number of networks requiring training and lowering the overall computational cost.
\end{remark}
 
The main steps of the PPTNN method are summarized in Algorithm \ref{alg:Parallel PhaseTNN}.

\begin{algorithm}[tb]
  \DontPrintSemicolon
  \SetKwInOut{Input}{Input}\SetKwInOut{Output}{Output}
  \SetKwProg{Fn}{Function}{}{}
  
  \Input{Training data $\{(x_{\xi},f(x_{\xi}))\}_{\xi=1}^{N_{\xi}}$, frequency nodes $\kappa_j = j \Delta k$, $j \in \mathcal{I}_{K}$, bandwidth $\Delta k$, threshold $\varepsilon$}
  \Output{The approximate solution $\tilde{f}(x)$}
  \BlankLine

  Initialize $\tilde{f}(x) \gets 0$\;
  \For{each $j \in \mathcal{I}_{K}$}{

    Compute the frequency-shifted data $\tilde{f}_j^{\text{shift}}(x_{\xi})$ using numerical quadrature according to Eq. \eqref{eq:shfittraningdate}\;
    Decompose $\tilde{f}_j^{\text{shift}}(x_{\xi})$ into real and imaginary parts: $\tilde{f}_j^{\text{real}}(x_{\xi})$, $\tilde{f}_j^{\text{imag}}(x_{\xi})$\;
    Compute root mean square (RMS) for each component: $\text{RMS}(\tilde{f}_j^{\text{real}})$, $\text{RMS}(\tilde{f}_j^{\text{imag}})$\;
    \BlankLine

    \For{$\vartheta \in \{\text{real}, \text{imag}\}$}{
      \If{$\text{RMS}(\tilde{f}_j^{\vartheta}) < \varepsilon$}{
        $\tilde{f}_j^{\vartheta}(x) \gets 0$\;
      }
      \Else{
        Train TransNet $T_j^{\vartheta}(x, \boldsymbol{\alpha}_j) $ to approximate $\tilde{f}_j^{\vartheta}(x)$\;
        $\tilde{f}_j^{\vartheta}(x) \gets T_j^{\vartheta}(x, \boldsymbol{\alpha}_j)$\;
      }
    }
    $\tilde{f}(x) \gets \tilde{f}(x) + [\tilde{f}_j^{\text{real}}(x) + i \tilde{f}_j^{\text{imag}}(x)] e^{2\pi i \kappa_j x}$\;
  }
  \Return $\tilde{f}(x)$\;
  
  \caption{PPTNN training scheme for 1D function approximation.}
  \label{alg:Parallel PhaseTNN}
\end{algorithm}

\subsection{A Parallel Phase-Shift Transferable Neural Network for the 2D Case}
In this subsection, we extend the PPTNN framework to the two-dimensional (2D) setting, generalizing the 1D high-frequency function approximation method to the domain $\Omega=[-1,1]^2$.

The PPTNN framework extends naturally to two spatial dimensions. For the spatial variable 
\(\boldsymbol{x}=(x_1,x_2) \in \mathbb{R}^2\) and frequency variable \(\boldsymbol{k}=(k_1,k_2)\in \mathbb{R}^2\), the 2D Fourier transform and its inverse are defined as
\begin{align}
\mathcal{F}[f](\boldsymbol{k}) = \int_{-\infty}^{\infty}\int_{-\infty}^{\infty} f(x_1, x_2)e^{-2\pi i (k_1 x_1 + k_2 x_2)}\,d x_1 d x_2,\\
\mathcal{F}^{-1}[\hat f](\boldsymbol{x}) = \int_{-\infty}^{\infty}\int_{-\infty}^{\infty} \hat f(k_1,k_2)e^{2\pi i (k_1 x_1 + k_2 x_2)}\,d k_1 d k_2.
\end{align}
We introduce a 2D frequency grid $(\kappa_{j_1},\kappa_{j_2})$, where
\begin{align}
\kappa_{j_1} = j_1 \Delta k_1,\quad \kappa_{j_2} = j_2 \Delta k_2,\qquad 
j_1 \in \mathcal{I}_{K_1}, \; j_2 \in \mathcal{I}_{K_2}. 
\end{align}
A partition of unity (POU) $\{\phi_{j_1,j_2}(k_1,k_2)\}$ for the rectangular frequency domain $[-K_1\Delta k_1,K_1\Delta k_1]\times[-K_2\Delta k_2,K_2\Delta k_2]$ is constructed as the tensor product of the 1D basis functions:
\begin{align}
\phi_{j_1,j_2}(\boldsymbol{k})=\phi\Bigl(\frac{k_1-\kappa_{j_1}}{\Delta k_1}\Bigr) \phi\Bigl(\frac{k_2-\kappa_{j_2}}{\Delta k_2}\Bigr),~~\text{and}\qquad 
\sum_{j_1,j_2}\phi_{j_1,j_2}(\boldsymbol{k})=1.
\end{align}
The inverse Fourier transform of
\(\phi_{j_1,j_2}(k_1,k_2)\) yields the 2D filtering kernel:
\begin{align}
\phi_{j_1,j_2}^\vee(\boldsymbol{x})=\Delta k_1 \Delta k_2 \phi^\vee\Bigl(\Delta k_1 x_1\Bigr)\phi^\vee\Bigl(\Delta k_2 x_2\Bigr)e^{2\pi i(\kappa_{j_1} x_1+\kappa_{j_2} x_2)}.
\end{align}
Using this kernel, we decompose \(f\) into frequency components via convolution:
\begin{align}
f_{j_1, j_2}(\boldsymbol{x}) = (\phi_{j_1, j_2}^{\vee} * f)(\boldsymbol{x}), \qquad 
f(\boldsymbol{x}) = \sum_{j_1=-K_1}^{K_1} \sum_{j_2=-K_2}^{K_2} f_{j_1, j_2}(\boldsymbol{x}). 
\label{eq:2D f_j}
\end{align}
Each component \(f_{j_1, j_2}\) has its spectrum confined to the rectangle \([\kappa_{j_1} - \frac{\Delta k_1}{2}, \kappa_{j_1} + \frac{\Delta k_1}{2}] \times [\kappa_{j_2} - \frac{\Delta k_2}{2}, \kappa_{j_2} + \frac{\Delta k_2}{2}]\).  
A phase-shift operation translates this spectrum to the baseband interval \([-\frac{\Delta k_1}{2}, \frac{\Delta k_1}{2}] \times [-\frac{\Delta k_2}{2}, \frac{\Delta k_2}{2}]\):
\begin{align}
f_{j_1, j_2}^{\text{shift}}(\boldsymbol{x}) = e^{-2\pi i (\kappa_{j_1} x_1 + \kappa_{j_2} x_2)} f_{j_1, j_2}(\boldsymbol{x}). 
\end{align}
To compute the training data for approximating  $f_{j_1,j_2}^{\text{shift}}(\boldsymbol{x}_{\xi})$, we introduce an auxiliary function $F(\boldsymbol{x})$ that ensures rapid spectral decay outside $\Omega$:
\begin{equation}
F(\boldsymbol{x}) = f(\boldsymbol{x}) \cdot \exp\left(-10\left[\left(\max(0,|x_1|-1)\right)^{6} + \left(\max(0,|x_2|-1)\right)^{6}\right]\right).
\label{Fzf2D}
\end{equation}

The training dataset $\{\boldsymbol{x}_{\xi},f_{j_1,j_2}^{\text{shift}}(\boldsymbol{x}_{\xi})\}_{\xi=1}^{N_{\xi}}$ 
for the 2D scenario is obtained via numerical quadrature:
\begin{equation}
\begin{aligned}
f_{j_1, j_2}^{\text{shift}}(\boldsymbol{x}_\xi) &\approx e^{-2\pi i \boldsymbol{\kappa}_{j_1, j_2}^\top \boldsymbol{x}_\xi} \sum_{s_1=1}^{N_{s_1}} \sum_{s_2=1}^{N_{s_2}} \phi_{j_1, j_2}^\vee(\boldsymbol{x}_\xi - \boldsymbol{x}_{s_1, s_2}) F(\boldsymbol{x}_{s_1, s_2}) \Delta\tilde{x}_{1, s_1} \Delta\tilde{x}_{2, s_2} \\
&:= \tilde{f}_{j_1, j_2}^{\text{shift}}(\boldsymbol{x}_\xi), \quad \xi \in \mathcal{I}_{N_{\xi}}^+.
\label{eq:shfittraningdate2D}
\end{aligned}
\end{equation}
Here, the grid spacings are defined as $\Delta x_1 = \frac{2C}{N_{s_1} - 1}$ and $\Delta x_2 = \frac{2C}{N_{s_2} - 1}$, respectively. The terms $\Delta \tilde{x}_{1, s_1}$ and $\Delta \tilde{x}_{2, s_2}$ denote the integration steps incorporating the trapezoidal rule weights, following the same formulation specified in Eq. \eqref{eq:trapz_weights}.

Analogous to the 1D case, we employ the Transferable Neural Network framework $T_{j_1,j_2}(\boldsymbol{x}, \boldsymbol{\alpha}_j)$ to approximate each component $f_{j_1,j_2}^{\text{shift}}$ by minimizing the following least-squares loss function:
\begin{equation}
L_{j_1,j_2}(\boldsymbol{\alpha}_{j_1,j_2}) = \sum_{\xi=1}^{N_{\xi}} \left( f_{j_1,j_2}^{\text{shift}}(\boldsymbol{x}_{\xi}) - T_{j_1,j_2}(\boldsymbol{x}_{\xi}, \boldsymbol{\alpha}_{j_1,j_2}) \right)^2,j_1 \in \mathcal{I}_{K_1}, \; j_2 \in \mathcal{I}_{K_2}.
\label{determT2D}
\end{equation}
From the phase-shift identity \eqref{Fzf2D} and the network approximation $T_{j_1,j_2}$, we reconstruct the component approximation:
$$
 \tilde{f}_{j_1,j_2}(\boldsymbol{x})=T_{j_1, j_2}(\boldsymbol{x})
e^{2\pi i (\kappa_{j_1} x_1 + \kappa_{j_2} x_2)}.
$$
Finally, the overall approximation $\tilde{f}(\boldsymbol{x})$ of $f(\boldsymbol{x})$ on $\Omega$ is obtained by summing the component approximations:
\begin{equation}
\tilde{f}(\boldsymbol{x}) = \sum_{j_1=-K_1}^{K_1} \sum_{j_2=-K_2}^{K_2} \tilde{f}_{j_1,j_2}(\boldsymbol{x}),\quad \boldsymbol{x}\in \Omega.
\end{equation}

\subsection{A Coupled Phase-Shift Transferable Neural Network (CPTNN)}
\label{sec:CPTNN-approximation} 
While the PPTNN framework decomposes data into frequency components for parallelizable TransNet training, it relies on explicit convolution (cf. Eq. \eqref{eq:shfittraningdate}) to construct task-specific training data—a bottleneck that fundamentally limits its scalability to high-dimensional problems and large-scale datasets. To address this limitation, we introduce a Coupled Phase-Shift Transferable Neural Network,  which integrates frequency modulation directly into the neural network architecture. Rather than explicitly decomposing the training data, the CPTNN constructs a set of phase-modulated neural basis functions, enabling the simultaneous representation of multiple frequency components within a single model.

Specifically, for a one-dimensional function
$f$, we decompose $f$ as 
$$f(x) = \sum_{j=-K}^K  f_{j}(x),$$
where the spectrum of $f_j(x)$ is confined to the interval $\left[\kappa_j-\frac{\Delta k}{2}, \kappa_j +\frac{\Delta k}{2}\right]$. To convert $f_j(x)$ into a tractable low-frequency problem, we introduce the phase-shifted variable
\begin{align}
g_{j}(x) = e^{-2 \pi i \kappa_{j} x} f_{j}(x) ,
\end{align}
whose spectrum is now centered at zero, corresponding to a low-frequency problem. Inverting this transformation yields
\begin{align}
f_{j}(x) = e^{2\pi i \kappa_{j} x} g_{j}(x).
\end{align}
We then employ TransNet to approximate these low-frequency components $g_j$. Substituting these approximations back into the expression for $f$, we obtain the final approximation
\begin{equation}
f_{T}(x)=\sum_{j=-K}^K e^{2\pi i \kappa_{j} x} T_{j}(x).
\label{eq:cpTNN41d}
\end{equation}
Equation \eqref{eq:cpTNN41d}  admits a natural generalization to higher-dimensional functions. For simplicity, given a $d$-dimensional function $f$, its approximation $f_T$ takes the following mathematical form:
\begin{equation}
\label{eq:coupleTNN1}
f_T(\boldsymbol{x}) = \sum_{\zeta=1}^{N_{\kappa}} e^{2 \pi i \boldsymbol{\kappa}_\zeta^\top \boldsymbol{x}}T_\zeta(\boldsymbol{x},\boldsymbol{\alpha}_\zeta),
\end{equation}
where \( \boldsymbol{x} \in \Omega \subset \mathbb{R}^d \) denotes the input variable and \( f : \mathbb{R}^d \to \mathbb{R} \) is the target function to be approximated. The vectors \( \boldsymbol{\kappa}_\zeta\) represent prescribed frequency modes generated on a uniformly spaced grid, i.e.,
$\boldsymbol{\kappa}_{\zeta}\in\{-K\Delta k,(-K+1)\Delta k,...,K\Delta k\}^{d}$. Here, \( T_\zeta(\boldsymbol{x},\boldsymbol{\alpha}_\zeta) \in \mathcal{P}_{\text{NN}}^L\) are compact complex-valued TransNet modules parameterized by coefficient vectors \( \boldsymbol{\alpha}_\zeta \).

For real-valued target functions, a real-valued Fourier series is mathematically equivalent to a complex-valued one. We may therefore restrict our attention to sine and cosine expansions, expressed as:
\begin{equation}
\label{eq:copRTNN}
f_{T}(\boldsymbol{x}) = \sum_{\zeta=1}^{N_{\kappa}} \left(T_\zeta^{\cos} \cos(2\pi\boldsymbol{\kappa}_\zeta^\top \boldsymbol{x}) + T_\zeta^{\sin} \sin(2\pi\boldsymbol{\kappa}_\zeta^\top \boldsymbol{x})\right),
\end{equation}
where \( T_\zeta^{\cos} , T_\zeta^{\sin} \in \mathcal{P}_{\text{NN}}^L \) are real-valued TransNets.

This formulation is equivalent to employing the following set of global basis functions:
\begin{align}\label{eq:basicf1}
\Psi_{\zeta,m}^{\cos}(\boldsymbol{x}) = \cos(2\pi \boldsymbol{\kappa}_\zeta^\top \boldsymbol{x}) \sigma(\boldsymbol{w}_{\zeta,m}^{\cos \top} \boldsymbol{x} + b_{\zeta,m}^{\cos}),\\\label{eq:basicf2}
\Psi_{\zeta,m}^{\sin} (\boldsymbol{x}) = \sin(2\pi \boldsymbol{\kappa}_\zeta^\top \boldsymbol{x}) \sigma(\boldsymbol{w}_{\zeta,m}^{\sin \top} \boldsymbol{x} + b_{\zeta,m}^{\sin}),
\end{align}
for $\zeta = 1, \dots, N_{\kappa}$ and $m = 1, \dots, M_{\text{sub}}$. Consequently, the CPTNN approximation admits a linear expansion over these fixed global basis functions:
\begin{align}
f_T(\boldsymbol{x}) = \sum_{\zeta=1}^{N_{\kappa}}\sum_{m=1}^{M_{\text{sub}}} \left( \alpha_{\zeta,m}^{\cos} \Psi_{\zeta,m}^{\cos}(\boldsymbol{x}) + \alpha_{\zeta,m}^{\sin} \Psi_{\zeta,m}^{\sin}(\boldsymbol{x}) \right).
\label{eq:ubasisEx}
\end{align}
The weights and biases of $\Psi_{\zeta,m}^{\cos}(\boldsymbol{x})$ and $\Psi_{\zeta,m}^{\sin}(\boldsymbol{x})$ are generated using the same sampling scheme as TransNet. Let $\boldsymbol{w}_{\eta}$ and $b_{\eta}$ denote any weight $\boldsymbol{w}_{\zeta,m}^{\cos/\sin}$ and bias $b_{\zeta,m}^{\cos/\sin}$, respectively. These parameters are decomposed as follows:
\begin{align}
\begin{cases}
\boldsymbol{w}_{\eta} = \gamma_{\eta} \boldsymbol{a}_{\eta},\\
 b{\eta} = \gamma_{\eta} r_{\eta},
\end{cases}
\end{align}
with the specific sampling procedure given by:
\begin{equation}
\boldsymbol{a}_{\eta} = \frac{X_{\eta}}{\| X_{\eta} \|_2}, \quad \text{and} \quad r_{\eta} = U_{\eta}, \ {\eta} = 1, \dots, 2N_{\kappa}M_{\text{sub}}.
\label{eq:a_eta&r_eta}
\end{equation}
Here, $\boldsymbol{X}_{\eta}$ are i.i.d. standard Gaussian distribution, and $U_{\eta}$ are i.i.d. uniform random variables on $[0, 1]$. Similar to PPTNN, the shape parameter $\gamma_{\eta}$ is fixed at 2, eliminating the need for manual tuning.

Since all basis functions are fixed, model training for approximating $f$ reduces to a linear least-squares problem. The optimal coefficient vector $\boldsymbol{\alpha}^*$ is obtained by minimizing the mean-squared error loss function:
\begin{equation}
\label{eq:loss_numerical}
L(\boldsymbol{\alpha}) = \sum_{\xi=1}^{N_{s}} \left| f(\boldsymbol{x}_{\xi}) - f_T(\boldsymbol{x}_{\xi}) \right|^2,
\end{equation}
where $\{\boldsymbol{x}_{\xi}\}_{\xi=1}^{N_s}$ are $N_s$ sampling nodes in $\Omega$. 

\section{Solving Differential Equations with CPTNN}
\label{sec:methodology solving} 
In this section, we apply the CPTNN framework to solve the following boundary value problem (BVP): 
\begin{equation}
\begin{aligned}
\mathcal{L}(u(\boldsymbol{x})) &= f(\boldsymbol{x}), \quad \boldsymbol{x} \in \Omega, \\
\mathcal{B}(u(\boldsymbol{x})) &= \mathcal{C}(\boldsymbol{x}), \quad \boldsymbol{x} \in \partial\Omega,
\end{aligned}
\label{eq:PDE}
\end{equation}
where $\Omega \subset \mathbb{R}^{d}$ is a bounded domain, $\mathcal{L}$ is a differential operator, $\mathcal{B}$ is a linear boundary operator, $f(\boldsymbol{x})$ is the source function, and $\mathcal{C}(\boldsymbol{x})$ denotes the prescribed boundary condition.

As established in Eq. \eqref{eq:ubasisEx} in Section \ref{sec:CPTNN-approximation}, the approximate solution $u_T$ to the exact solution $u$ admits the simple representation
\begin{equation}
\label{eq:coupleTNNappu}
u_{T}(\boldsymbol{x})=\sum_{\eta=1}^{M_{T}}\alpha_{\eta}\Psi_{\eta}(\boldsymbol{x}).
\end{equation}
where \(\{\Psi_{\eta}\}_{\eta=1}^{M_T}\) are the phase-modulated neural basis functions defined in Eqs. \eqref{eq:basicf1}-\eqref{eq:ubasisEx}, and \(\boldsymbol{\alpha}=(\alpha_1,\alpha_2,\dots,\alpha_{M_T})^\top\)
is the vector of unknown coefficients to be determined.

The network parameters $\boldsymbol{\alpha}$ are optimized by minimizing a composite loss function that enforces satisfaction of the PDE and boundary conditions at collocation points:
\begin{equation}
\begin{aligned}
L(\boldsymbol{\alpha}) = \sum_{\mu=1}^{N_{\text{pde}}} \left| f(\boldsymbol{x}_{\mu}) - \mathcal{L}(u_T(\boldsymbol{x}_{\mu})) \right|^2 + \sum_{{\nu}=1}^{N_{\text{bc}}} \left| \mathcal{C}(\boldsymbol{y}_{\nu}) - \mathcal{B}(u_T(\boldsymbol{y}_{\nu})) \right|^2. 
\end{aligned}
\end{equation}
where $\{\boldsymbol{x}_{\mu}\}_{{\mu}=1}^{N_{\text{pde}}}\subset \Omega$ and $\{\boldsymbol{y}_{\nu}\}_{{\nu}=1}^{N_{\text{bc}}}\subset\partial\Omega$ denote the sets of interior and boundary collocation points, respectively.
\paragraph{Linear Differential Operator}
When $\mathcal{L}$ is a linear differential operator, the linearity of $\mathcal{L}$ and $\mathcal{B}$ implies 
\begin{align}
\mathcal{L}(u_T(\boldsymbol{x})) = \sum_{\eta=1}^{M_{T}} \alpha_{\eta} \mathcal{L}[\Psi_{\eta}(\boldsymbol{x})].
\end{align}
We define the feature matrices $\mathbf{A}_{\text{pde}} \in \mathbb{R}^{N_{\text{pde}} \times M_{T}}$ and $\mathbf{A}_{\text{bc}} \in \mathbb{R}^{N_{\text{bc}} \times M_{T}}$  with entries
\begin{align}
(\mathbf{A}_{\text{pde}})_{\mu,\eta} = \mathcal{L}[\Psi_{\eta}(\boldsymbol{x}_{\mu})], \quad (\mathbf{A}_{\text{bc}})_{\nu,\eta} = \mathcal{B}[\Psi_{\eta}(\boldsymbol{y}_{\nu})]. 
\end{align}
Let $\mathbf{f} \in \mathbb{R}^{N_{\text{pde}} \times 1}$ and $\mathbf{c} \in \mathbb{R}^{N_{\text{bc}} \times 1}$ denote the target vectors with entries $(\mathbf{f})_{\mu} = f(\boldsymbol{x}_{\mu})$ and $(\mathbf{c})_{\nu} = \mathcal{C}(\boldsymbol{y}_{\nu})$, respectively.
The optimal coefficient vector $\boldsymbol{\alpha}$ is then obtained by solving the linear least-squares problem:
\begin{align}
\bm{\alpha}^{\text{OPT}} = \arg\min_{\bm{\alpha}} \left\| \begin{bmatrix} \mathbf{A}_{\text{pde}} \\ \mathbf{A}_{\text{bc}} \end{bmatrix} \bm{\alpha} - \begin{bmatrix} \mathbf{f} \\ \mathbf{c} \end{bmatrix} \right\|_2^2.
\end{align}

\paragraph{Nonlinear Differential Operator}
For nonlinear $\mathcal{L}$, we employ a \textbf{Picard iterative method} to linearize the problem. We decompose $\mathcal{L}$ into a linear component $\mathcal{A}$ and a nonlinear component $\mathcal{N}$, i.e., $\mathcal{L}(u) = \mathcal{A}(u) + \mathcal{N}(u)$. At the $\text{iter}$-th iteration, the nonlinear term $\mathcal{N}(u)$ is evaluated using the previous iterate $u^{(\text{iter})}(\boldsymbol{x})$, reducing the problem to a sequence of linear least-squares subproblems. The updated coefficients $\bm{\alpha}^{(\text{iter}+1)}$ are computed by minimizing
\begin{align}
\bm{\alpha}^{(\text{iter}+1)} = \arg\min_{\bm{\alpha}} \left\| \begin{bmatrix} \mathbf{A}_{\text{pde}} \\ \mathbf{A}_{\text{bc}} \end{bmatrix} \bm{\alpha} - \begin{bmatrix} \mathbf{f} - \mathbf{N}(\bm{\alpha}^{(\text{iter})}) \\ \mathbf{c} \end{bmatrix} \right\|_2^2,
\end{align}
where $\mathbf{A}_{\text{pde}}$ and $\mathbf{A}_{\text{bc}}$ are the feature matrices associated with $\mathcal{A}$ and $\mathcal{B}$, and $\mathbf{N}(\boldsymbol{\alpha}^{(\text{iter})}) \in \mathbb{R}^{N_{\text{pde}} \times 1}$ is the discrete nonlinear vector with entries $(\mathbf{N}(\boldsymbol{\alpha}^{(\text{iter})}))_{\mu} = \mathcal{N}(u^{(\text{iter})}(\boldsymbol{x}_{\mu})).$ 

Iterations proceed until the relative change in the solution falls below a tolerance tol: 
\begin{align}
\frac{\|u^{(\text{iter}+1)}-u^{(\text{iter})}\|_2}{\|u^{(\text{iter})}\|_2} < \text{tol}. 
\end{align}

\begin{remark}
The CPTNN framework also exhibits strong approximation performance for other classes of partial differential equations (PDEs) with high-frequency solutions, including the large-wavenumber Helmholtz equation and the time-dependent wave equation. We validate the successful application of our method to such PDEs in the numerical examples section.
\end{remark}
\begin{remark}
Notably, the proposed method admits a natural generalization to the approximation of discontinuous functions and the solution of interface problems. Following the domain-decomposition-based TransNet framework established in the literature, we construct separate neural networks on individual subdomains to address these challenging problems. To rigorously validate this generalization, we test the method on the piecewise function $f_3$ in our numerical experiments, and the results confirm that our approach maintains high-precision approximation even for discontinuous target functions.
\end{remark}
\begin{remark}
Applying PPTNN to solve the boundary value problem (4.1) suffers from essential and insurmountable difficulties. The key bottleneck is that, in decomposing the original equation into frequency-wise subproblems, the boundary conditions of the subproblems lack a unified, rigorous construction and are extremely difficult to properly determine.
\end{remark}

\section{Numerical Results}
\label{sec:numerical results}

In this section, we present numerical results to demonstrate the superior performance of the proposed methods. Our examples encompass both function approximation and partial differential equation (PDE) solving. We employ the relative $L_{2}$ error to quantify performance, defined as
$$\frac{\left\|u_{\theta}-u\right\|_{2}}{\left\|u\right\|_{2}},$$
where $u$ denotes the exact solution and $u_{\theta}$ the corresponding approximate solution.

All neural networks in our experiments utilize the tanh activation function. All least-squares problems are solved using the  \texttt{numpy.linalg.lstsq} routine. Computational timing data is collected on a local workstation (2.90 GHz Intel Core Ultra 9 285H CPU, 32GB memory).

For the Transferable Neural Network, the shape parameter tuning strategy proposed in \cite{zhangTransferableNeuralNetworks2024} is no longer applicable. We instead identify the optimal shape parameter via a systematic grid search, selecting the approximation with the minimal error as the baseline for comparison with our proposed method.

For the random feature method (RFM) \cite{chenBridgingTraditionalMachine2022}, hidden-layer network coefficients are initialized with uniform random values sampled from the interva $[-R_{\max}, R_{\max}]$, where $R_{\max}$ is tuned for each specific high-frequency problem. Analogous to the TransNet setup, we determine the optimal hyperparameters via grid search and use the resulting minimal-error approximation as the reference baseline for evaluating our method.

In the following numerical experiments, the number of hidden neurons for both RFM and TransNet is set to an exceptionally large value of 100,000, unless otherwise specified. This over-parameterization is adopted to strictly eliminate potential performance limitations arising from insufficient network capacity, ensuring a fair comparison across all methods.

\subsection{Filtering Accuracy}
\label{sec:filtering accuracy}

In this subsection, we compare the accuracy of the filtering method from Ref. \cite{caiPhaseShiftDeep2020} with the two novel filtering strategies proposed in this work. The experimental setup is detailed as follows:
 
\begin{itemize}
    \item \textbf{Method 1.} The filtering method from \cite{caiPhaseShiftDeep2020}, applied directly to the raw data. 

    \item \textbf{Method 2.} Based on the raw data, an auxiliary function is first constructed using Eq. \eqref{eq:hf1}, followed by filtering via Eq. \eqref{eq:traningdate}.  

    \item \textbf{Method 3.} Based on the raw data, an auxiliary function is first constructed using Eq. \eqref{eq:hf2}, followed by filtering via Eq. \eqref{eq:traningdate}.
\end{itemize}

\begin{table}[!h]
\centering
\caption{Comparison of the performance of three filtering methods.}
\label{tab:filtering results}
\begin{tabular}{c|c|c|c}
\hline
Function & Method 1 & Method 2 & Method 3 \\
\hline
$x^{3}$ & 6.31E-03 & 8.19E-15 & 8.17E-15 \\
\hline
$100^x$ & 3.52E+04 & 2.16E-14 & 2.11E-14 \\
\hline
$\sin(2\pi x)$ & 3.42E-06 & 1.79E-15 & 2.69E-15 \\
\hline
$\sin(2x + 1) + 0.2e^{1.3x}$ & 1.17E-02 & 2.45E-15 & 2.45E-15 \\
\hline
$(1-0.5x^{2})\cos(30(x+0.5x^{3}))$ & 6.40E-02 & 3.05E-15 & 5.96E-08 \\
\hline
$x\sin(100x)$ & 9.73E-04 & 1.09E-14 & 1.04E-06 \\
\hline
\end{tabular}
\end{table}
We report the filtering errors in Table \ref{tab:filtering results}, where the \textbf{filtering error} norm: 
\begin{equation}
\textbf{filtering error} = \frac{\left\| f(x_{\xi}) -\sum_{j=-K}^{K}\tilde{f}_j(x_{\xi}) \right\|_2}{\left\| f(x_{\xi}) \right\|_2},
\label{eq:errorfilter}
\end{equation}
with $K=40$ (corresponding to $C=5$ in Eq. \eqref{eq:traningdate} for computing $\tilde{f}_j$). All methods are evaluated using 5001 sampling points. The bandwidth
$\Delta k$ is fixed at 2, and the error in Eq. \eqref{eq:errorfilter} is evaluated over 1001 equidistant points $\{x_{\xi}\}_1^{1001}\subset[-1,1]$. 

Numerical results in Table \ref{tab:filtering results} demonstrate that the two proposed filtering strategies (Methods 2 and 3) significantly enhance the filtering performance compared to Method 1 from \cite{caiPhaseShiftDeep2020}, with both achieving accuracy close to machine precision. This study provides a necessary foundation for the subsequent application of the PPTNN method to high-precision function approximation and PDE solving.

\subsection{Accuracy of approximation of functions with PhaseTNN}
In this subsection, we perform numerical validation to evaluate the approximation performance of the proposed PhaseTNN framework (encompassing both PPTNN and CPTNN) for high-frequency functions. We begin by considering the following three one-dimensional test functions, which are adopted from \cite{wangMultistageNeuralNetworks2024}, \cite{xuMultigradeDeepLearning2025}, and \cite{caiPhaseShiftDeep2020}, respectively:

\begin{align*}
f_1(x) &= (1-0.5x^{2})\cos\bigl(a(x+0.5x^{3})\bigr),\\
f_2(x) &= x\sin(100x),\\
f_3(x) &= \begin{cases} 
\sin(x) + \sin(3x), & x \in [-1, 0), \\ 
\sin(23x) + \sin(137x) + \sin(203x), & x \in [0, 1].
\end{cases}
\end{align*}

Subsequently, we extend our evaluation to two-dimensional high-frequency functions, using the following test cases:
\begin{align*}
f_{4}(x_1, x_2) &= \left(1 - 0.5 x_1^2\right) \cos\bigl[20(x_1 + 0.5 x_1^3)\bigr] \cdot (1 - 0.5 x_2^2) \cos\bigl[20(x_2 + 0.5 x_2^3)\bigr] , \\
f_{5}(x_1, x_2) &= x_1 \cos\bigl(2 \pi \cdot 5 (x_1 + x_2)\bigr)+ (x_1 x_2)^3 + e^{\sin(7 x_1 x_2)}. 
\end{align*}

For the one–dimensional problems, we adopt 1001 equidistant points in \([-1, 1]\) as training data, and evaluate approximation accuracy via the relative 
\(L_2\) error computed on 8000 test points.
For the 2D function approximation, we set the domain to $[-1, 1]^2$, with training and test data consist of $201^2$ and $300^2$ equidistant points, respectively.

\subsubsection{Effects of the shape parameter and network capacity}
The frequency content of $f_1$ increases monotonically with $a$. To assess performance across frequency regimes, we select \( a = 1, 30, 100 \) (low, moderate, high) and conduct approximation experiments with TransNet, PPTNN, and CPTNN. Network configurations are as follows: TransNet and CPTNN employ 1000 hidden neurons, while PPTNN uses 100 neurons per sub-network.

\begin{figure}[htbp]
    \centering
    \includegraphics[width=0.9\textwidth]{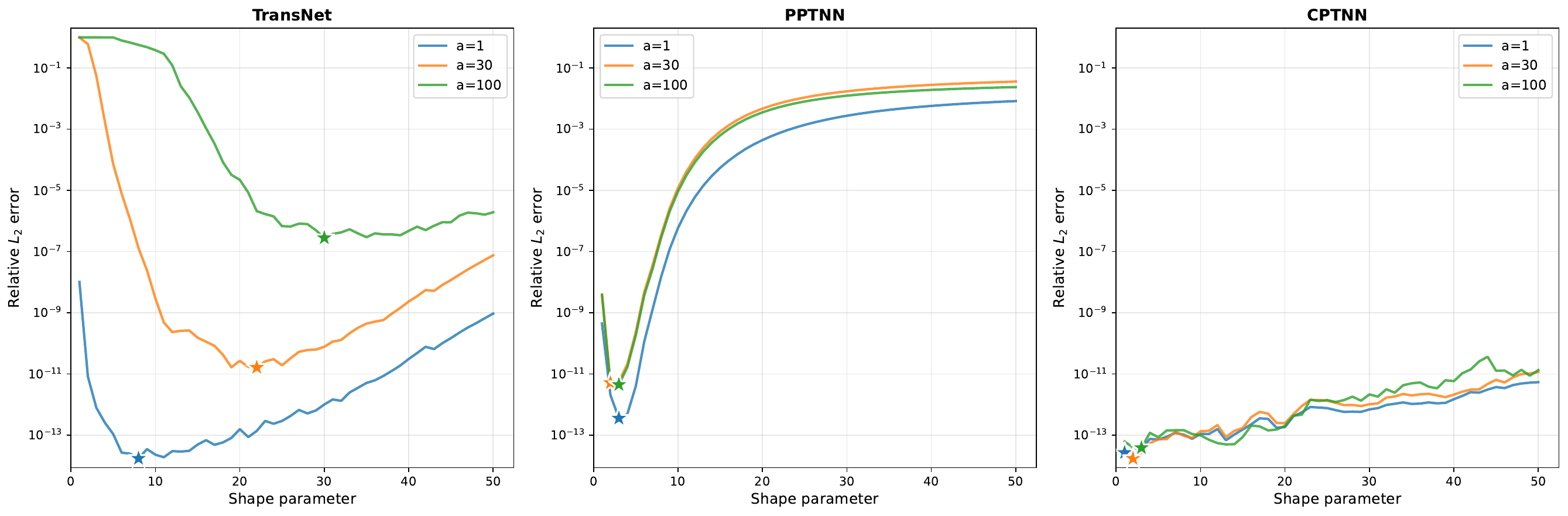}
    \caption{Relative $L_2$ error vs. $\gamma$ for $f_1$ (\( a = 1, 30, 100 \)) with TransNet, PPTNN and CPTNN.}
    \label{fig:1D M1&30&100 TNN&PPTNN&CPTNN Error&Shape}
\end{figure}

We evaluate the sensitivity  of approximation accuracy to $\gamma$. Fig. \ref{fig:1D M1&30&100 TNN&PPTNN&CPTNN Error&Shape} plots relative \( L_2 \) test error  error vs. $\gamma$, with star markers denoting optimal. For TransNet, the optimal  $\gamma$ fluctuates sharply with frequency, requiring costly parameter tuning. In contrast, PPTNN and CPTNN maintain stable optimal  $\gamma$, streamlining implementation and improving robustness.

\begin{figure}[htbp]
    \centering
    \includegraphics[width=0.9\textwidth]{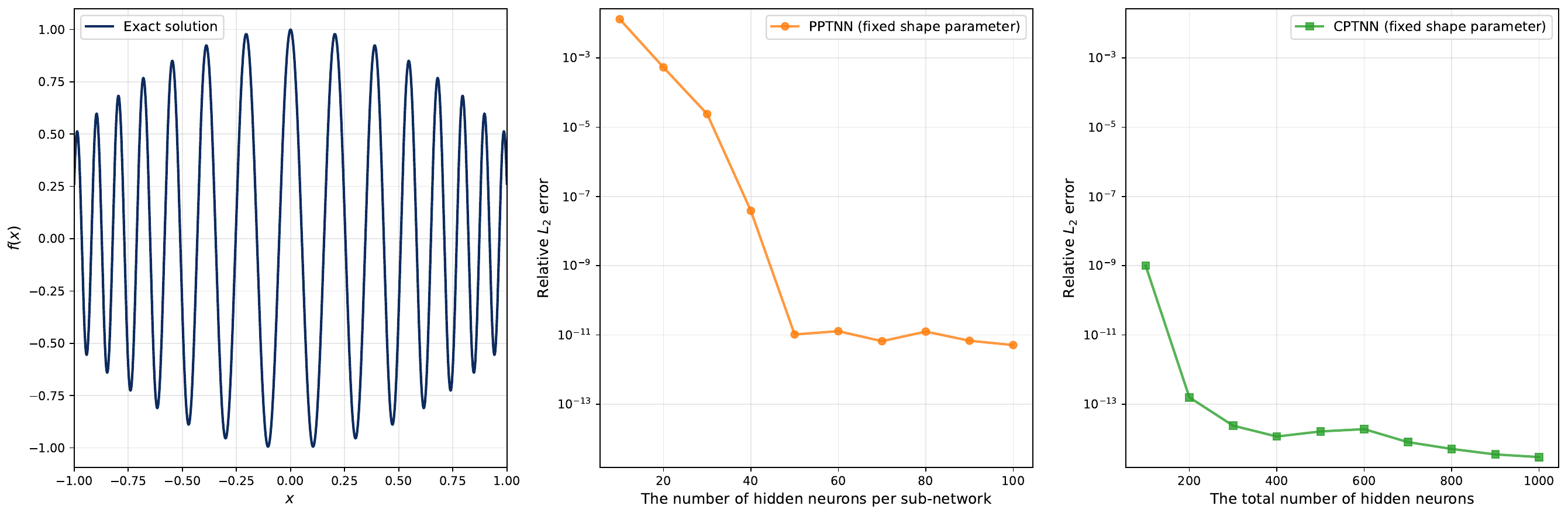}
    \caption{Convergence behavior for $f_1$ ($a=30$). (Left) $f_1$; (Middle) PPTNN relative $L_2$ error vs. hidden neurons per sub-network; (Right) CPTNN relative $L_2$ error vs. total hidden neurons.}
    \label{fig:1D 1D F1-A30 PPTNN&CPTNN Error&Basis}
\end{figure}

We further analyze convergence behavior for the highly oscillatory function $f_1(x)$ ($a=30$) with $\gamma=2$ fixed. PPTNN uses $N_s=5001$, $\Delta k=2$, $\kappa_j \in \{-50, -48, \dots, 50\}$, $\varepsilon=10^{-14}$, varying sub-network neurons (10-100). CPTNN varies total hidden neurons (100-1000), corresponding to 25-250 frequency groups. As shown in Fig. \ref{fig:1D 1D F1-A30 PPTNN&CPTNN Error&Basis}, both methods exhibit monotonic error decay with increasing capacity, demonstrating strong expressivity after low-frequency decomposition. Notably, PPTNN’s error saturates at
$\mathcal{O}(10^{-12})$ due to summation-induced round-off errors in its parallel architecture. Conversely,  CPTNN resolves all frequency components in a coupled manner, avoiding such error accumulation and achieving higher precision ($\mathcal{O}(10^{-15})$).

\subsubsection{1D function approximation}
We compare our proposed methods against four existing algorithms--FNN \cite{hornikMultilayerFeedforwardNetworks1989}, RFM \cite{chenBridgingTraditionalMachine2022}, MSNN \cite{wangMultistageNeuralNetworks2024} and TransNet \cite{zhangTransferableNeuralNetworks2024}--for the one--dimensional problems $f_1$, $f_2$  and $f_3$. 

Training details for $f_1$ and $f_2$ are summarized in Table \ref{tab:1D configurations}. For the piecewise function $f_3$, we use domain decomposition, training two sub-networks on $[-1,0)$ and $[0,1]$ separately. This strategy is applied to all models for fairness, with Table \ref{tab:1D configurations} listing per-sub-network configurations. Minor adjustments are made for $f_3$: FNN architecture per sub-network is $[1,16100,1]$, and CPTNN uses 250 frequencies sampled from $[0,40]$.

\begin{table}[!h]
    \centering
    \caption{Network architectures and training strategies for the 1D function approximation.}
    \label{tab:1D configurations}
    \begin{tabular}{c|c|c|c}
        \hline
        Method & Architecture & Training strategy & Other configurations \\
        \hline
        FNN \cite{hornikMultilayerFeedforwardNetworks1989} & $[1, 10100, 1]$ & Adam: 10k (steps) & \makecell{Learning rate: 0.001 \\ (decay 0.7 per 2k steps)} \\
        \hline
        RFM \cite{chenBridgingTraditionalMachine2022} & $[1, 100000, 1]$ & Least squares & $R_{\max} \in \{2, 4, \dots, 200\}$ \\
        \hline
        MSNN \cite{wangMultistageNeuralNetworks2024} & \makecell{Stage 1-2: $[1, 20 \times 3, 1]$ \\ Stage 3-4: $[1, 30 \times 3, 1]$} & \makecell{Adam + L-BFGS \\ Stage 1: 3k + 10k (steps) \\ Stage 2: 5k + 20k (steps) \\ Stage 3: 5k + 30k (steps) \\ Stage 4: 5k + 40k (steps)} & \makecell{Stage 1-2: $N_{\text{train}} = 1501$ \\ Stage 3-4: $N_{\text{train}} = 4801$} \\
        \hline
        TransNet \cite{zhangTransferableNeuralNetworks2024} & $[1, 100000, 1]$ & Least squares & $\gamma \in \{2, 4, \dots, 100\}$ \\
        \hline
        PPTNN & \makecell{$[1, 100, 1]$ \\ (per sub-net)} & Least squares & \makecell{$N_s=5001$ \\ $C=5$ \\ $\Delta k=2$ \\ $\kappa_j \in [-50, 50]$ \\ $\varepsilon=10^{-14}$} \\
        \hline
        CPTNN & $[1, 1000, 1]$ & Least squares & 250 frequencies in $[0, 20]$ \\
        \hline
    \end{tabular}
\end{table}

Table \ref{tab:Approximation 1D f1&f2&f3 error&time} comprehensively evaluates the proposed frameworks against representative baselines, assessing both approximation accuracy and computational efficiency. As shown in the upper panel of Table \ref{tab:Approximation 1D f1&f2&f3 error&time}, traditional architectures such as FNN and RFM fail to converge for highly oscillatory functions $f_1$ and $f_2$, with errors stagnating at
$\mathcal{O}(10^0)$. While the multi-stage neural network (MSNN) partially mitigates this issue, it still suffers from a bottleneck in high-frequency approximation, saturating at  $\mathcal{O}(10^{-6})$--$\mathcal{O}(10^{-7})$. Under the optimal shape parameter, TransNet achieves strong performance across all benchmarks, attaining relative $L_2$ errors in the range $\mathcal{O}(10^{-13})$--$\mathcal{O}(10^{-11})$. PPTNN matches TransNet’s accuracy, whereas CPTNN consistently outperforms all compared methods, reducing errors to near-machine precision ($\mathcal{O}(10^{-14})$--$\mathcal{O}(10^{-15})$) for both high- and low-frequency functions.

When evaluating computational efficiency (lower panel of Table \ref{tab:Approximation 1D f1&f2&f3 error&time}), the practical cost of achieving high accuracy is critical. Existing high-accuracy methods incur prohibitive overheads: MSNN requires over 1200 seconds due to its sequential multi-stage training. For TransNet, optimal accuracy relies on an exhaustive grid search over the shape parameter $\gamma$ (2 to 50 in steps of 2, totaling 25 independent runs). Although this tuning could be parallelized, it wastes $96\%$ of computational effort. To accurately reflect the practical computational budget, the time reported for TransNet encompasses the total duration of shape parameter tuning (exceeding 56 seconds), whereas for PPTNN, it denotes the theoretical parallel execution time. This clearly demonstrates the resource inefficiency of parameter-dependent models.

In stark contrast, our proposed PPTNN and CPTNN resolve this accuracy–efficiency dilemma by eliminating exhaustive parameter searches. By decomposing the target function in the frequency domain, PPTNN enables fully independent training of its sub-networks. With sufficient parallel computing resources, its execution time is governed by the maximum runtime across sub-networks:
$$
t_{\text{parallel}} = \max(t_{\text{filter}}) + \max(t_{\text{train}}),
$$
yielding an unprecedented runtime of less than 0.6 seconds. Furthermore, CPTNN achieves the highest precision in approximately 1 second. By removing the need for hyperparameter tuning, our frameworks deliver speedups of several orders of magnitude over tuning-dependent baselines while retaining extreme precision.

\begin{table}[!h]
    \centering
    \caption{Comparison of approximation accuracy and computational efficiency for 1D functions}
    \label{tab:Approximation 1D f1&f2&f3 error&time}
    \begin{tabular}{c|c|c|c|c|c|c}
        \hline
        \multicolumn{7}{c}{Relative $L_2$ error} \\ 
        \hline
        Function & FNN \cite{hornikMultilayerFeedforwardNetworks1989} & RFM \cite{chenBridgingTraditionalMachine2022} & MSNN \cite{wangMultistageNeuralNetworks2024} & TransNet \cite{zhangTransferableNeuralNetworks2024} & PPTNN & CPTNN \\
        \hline
        $f_1$   & 9.99E-01 & 7.39E-11 & 8.29E-07 & 3.57E-11 & 5.13E-12 & 2.98E-15 \\
        $f_2$  & 1.00E+00 & 4.05E-10 & 3.74E-06 & 1.42E-10 & 8.12E-12 & 2.55E-14 \\
        $f_3$   & 7.01E-01 & 1.96E-09 & 1.28E-06 & 3.00E-09 & 7.41E-13 & 1.54E-13 \\
        \hline
        \multicolumn{7}{c}{Time (seconds)} \\
        \hline
        Function & FNN & RFM & MSNN & TransNet & \makecell{PPTNN \\ (parallel time)} & CPTNN \\
        \hline
        $f_1$   & 205.30 & 10.27 & 1288.37 & 11.01 & 0.36 & 0.99 \\
        $f_2$ & 190.76 & 12.93 & 1392.49 & 13.04 & 0.54 & 0.88 \\
        $f_3$    & 326.15 & 7.48 & 1963.08 & 7.55 & 0.21 & 0.78 \\
        \hline
    \end{tabular}
\end{table}

\subsubsection{2D function approximation}
To evaluate the 2D function approximation for $f_4$,  we fix the shape parameter at $\gamma=2$. For PPTNN, we set $N_s=1001^2$, $C=5$, $\Delta k=2$, $\varepsilon=10^{-12}$, and the frequency centers $\boldsymbol{\kappa}_j \in \{-40, -38, \dots, 40\}^2$, while varying the number of neurons per sub-network ($M_{\text{sub}}$) from 50 to 2000. For CPTNN, we set the frequency range to $[-10, 10]^2$, and vary the total number of hidden neurons ($M_T$) from 1600 ($4 \times 20^2$) to 10000 ($4 \times 50^2$).

Table \ref{tab:Approximation 2D f4 error&basis} summarizes the convergence behavior of PPTNN and CPTNN for the 2D oscillatory function $f_4$. As network capacity increases, both methods exhibit a decay in relative $L_2$ error,  with values falling to $\mathcal{O}(10^{-10})$. This confirms that the superior approximation capabilities observed in 1D experiments are preserved in 2D problems.

\begin{table}[!h]
    \centering
    \caption{Convergence behavior of PPTNN and CPTNN for the 2D function $f_4$.}
    \label{tab:Approximation 2D f4 error&basis}
    \begin{tabular}{c|c|c|c}
        \hline
        \multicolumn{2}{c|}{PPTNN} & \multicolumn{2}{c}{CPTNN} \\
        \hline
        $M_{\text{sub}}$ & Relative $L_{2}$ error & $M_T$ & Relative $L_{2}$ error \\
        \hline
        50   & 7.30E-02 & $4\times20^2$ & 3.01E-01 \\
        100  & 2.56E-03 & $4\times25^2$ & 1.14E-04 \\
        200  & 3.88E-05 & $4\times30^2$ & 9.80E-08 \\
        500  & 3.01E-08 & $4\times35^2$ & 8.02E-09 \\
        1000 & 2.11E-09 & $4\times40^2$ & 2.64E-09 \\
        1500 & 1.21E-09 & $4\times45^2$ & 1.25E-09 \\
        2000 & 9.87E-10 & $4\times50^2$ & 9.07E-10 \\
        \hline
    \end{tabular}
\end{table}

To further assess our approach, we compare the proposed PPTNN and CPTNN against FNN, RFM, and TransNet for the 2D high-frequency functions $f_4$ and $f_5$. Detailed network configurations are summarized in Table \ref{tab:2D configurations}. Consistent with the 1D experiments, the baseline models are allocated substantial parameter capacities to ensure a fair comparison. Table \ref{tab:Approximation 2D f4&f5 Error} reports the relative $L_2$ test errors, confirming that our methods successfully extend their high-precision approximation capabilities to 2D domains. Notably, PPTNN and CPTNN consistently outperform all other compared methods.

\begin{table}[!h]
    \centering
    \caption{Network architectures and training strategies for the 2D function approximation.}
    \label{tab:2D configurations}
    \begin{tabular}{c|c|c|c}
        \hline
        Method & Architecture & Training strategy & Other configurations \\
        \hline
        FNN & $[2, 100\times5, 1]$ & Adam (10000 steps) & \makecell{Learning rate: 0.001 \\ (decay 0.7 per 2k steps)} \\
        \hline
        RFM & $[2, 100000, 1]$ & Least squares & $R_{\max} \in \{2, 11, \dots, 101\}$ \\
        \hline
        TransNet & $[2, 100000, 1]$ & Least squares & $\gamma \in \{2, 4, \dots, 30\}$ \\
        \hline
        PPTNN & \makecell{$[2, 2000, 1]$ \\ (per sub-net)} & Least squares & \makecell{$N_s=1001^2$ \\ $C=5$ \\ $\Delta k=2$ \\ $\boldsymbol{\kappa}_j \in \{-40,-38,\dots,40\}^2$ \\ $\varepsilon=10^{-12}$} \\
        \hline
        CPTNN & $[2, 10000, 1]$ & Least squares & \makecell{2500 frequencies in \\ $[-10, 10]^2$} \\
        \hline
    \end{tabular}
\end{table}

\begin{table}[!h]
    \centering
    \caption{Performance comparison of approximation accuracy for the 2D functions.}
    \label{tab:Approximation 2D f4&f5 Error}
    \begin{tabular}{c|c|c|c|c|c}
        \hline
        function & FNN & RFM & TransNet & PPTNN & CPTNN \\
        \hline
        \(f_4\) & 1.00E+00 & 3.45E-07 & 7.66E-08 & 9.87E-10 & 9.07E-10 \\
        \(f_5\) & 9.20E-03 & 3.26E-08 & 1.01E-08 & 3.69E-09 & 3.92E-12 \\
        \hline
    \end{tabular}
\end{table}

\subsection{PhaseTNN for solving PDEs with high frequency solutions}
In this subsection, we evaluate the performance of our CPTNN for solving PDEs with high-frequency solutions under the following consistent settings:
\begin{itemize}
    \item \textbf{1D problems}: Trained on 2001  equidistant points, tested on 8000 points for relative $L_2$ error.
    \item \textbf{2D problems}: Trained on a \(100 \times 100\) grid, tested on a finer   \(200 \times 200\) grid.
    \item \textbf{Baseline PINN}:  Single-hidden-layer network for fair comparison. Optimized via Adam for 10,000 steps, initial learning rate of $0.001$ that is decayed by a factor of $0.7$ every $2000$ steps.
\end{itemize}

\subsubsection{Variable-coefficient elliptic equation}
We first consider the following one-dimensional equation:
\begin{equation}
\begin{aligned}
-\bigl((1+x^2)u'(x)\bigr)' &= f(x), \quad x\in[-1,1] \\
u(-1) = u(1) &= 0.
\end{aligned}
\label{eq:variable coefficient}
\end{equation}
We choose the source term $f(x)$ such that Eq. \eqref{eq:variable coefficient} has the following solution:
\begin{equation}
u(x) = \sin(2\pi x) + 0.1\sin(a \pi x^2),
\label{eq:variable coefficient solution}
\end{equation}
with \(a>0\) a parameter controlling the high-frequency content of the problem.

In this numerical study, we fix the number of hidden neurons to
$4000$ for both PINN and CPTNN. For CPTNN, these neurons are constructed using $2000$ distinct frequency modes, uniformly sampled from the interval  $[0, 300]$ across all three test configurations ($a=100, 150, 250$).
For the baseline methods, we perform targeted parameter tuning as follows:
\begin{itemize}
    \item For RFM, the scaling parameter $R_{\max}$ is tuned over the set $\{1,11,\dots,301\}$.
    \item For TransNet, the shape parameter is exhaustively tuned over the set $\{2,4,\dots,350\}$.
\end{itemize}

\begin{figure}[htbp]
    \centering
    \includegraphics[width=0.6\textwidth]{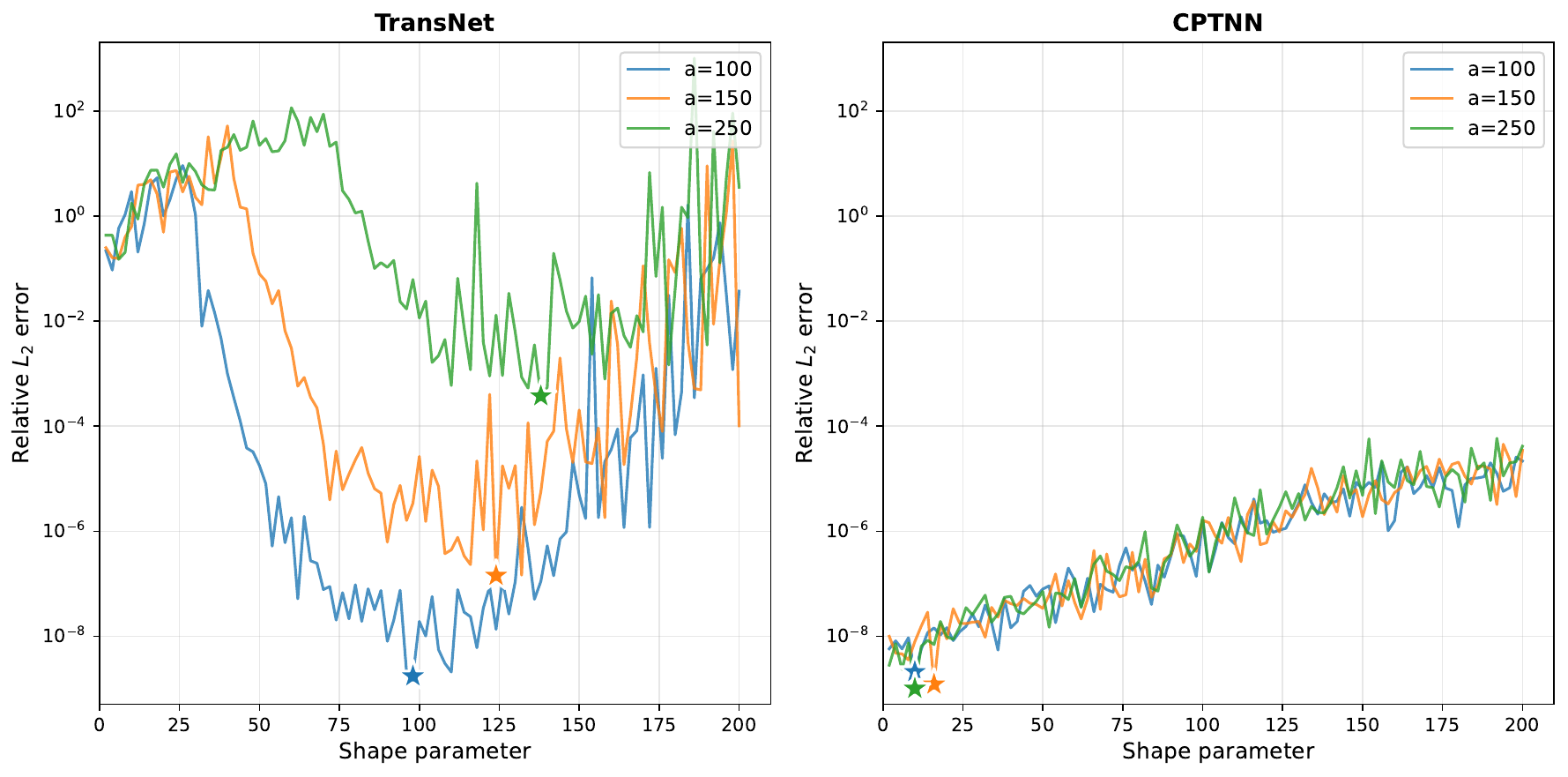}
    \caption{Relative $L_2$ error vs. shape parameter $\gamma$ for variable coefficient equation \eqref{eq:variable coefficient}}
    \label{fig:1D Variable coefficient TNN&CPTNN Error&Shape}
\end{figure}

To evaluate the influence of the shape parameter
$\gamma$ for TransNet and CPTNN, we vary  $\gamma$  over the set $\{2, 4, \dots, 200\}$ for both methods.
Figure \ref{fig:1D Variable coefficient TNN&CPTNN Error&Shape} depicts the relative $L_2$ error versus  $\gamma$ for TransNet and CPTNN under these frequency settings. Numerical results from Figure \ref{fig:1D Variable coefficient TNN&CPTNN Error&Shape}, demonstrate that the optimal shape parameter of CPTNN is nearly independent of the problem frequency, consistently  achieving high precision  $\mathcal{O}(10^{-9})$ within a narrow range of $\gamma$. This robust behavior rigorously justifies fixing the shape parameter to a small constant(e.g., $\gamma=2$): by intrinsically avoiding numerically ill-conditioned regions, CPTNN guarantees both high efficiency and high accuracy, regardless of the target function’s frequency. In contrast, the optimal shape parameter of TransNet exhibits a strong dependence on the problem frequency, and its peak achievable accuracy deteriorates with increasing frequency.

Table \ref{tab:Variable coefficient Error&time} compares the relative $L_2$ errors and computational times of numerical solutions obtained by various methods for solving Eq. \eqref{eq:variable coefficient} with a high-frequency solution ($a=250$). Numerical results demonstrate that PINN fails to effectively solve this high-frequency problem. The relative errors of RFM and TransNet are
$\mathcal{O}(10^{-5})$ and $\mathcal{O}(10^{-6})$, respectively, whereas the proposed CPTNN method attains a relative $L_2$ error of  $\mathcal{O}(10^{-9})$ representing an improvement of \textbf{3-4 orders of magnitude} in accuracy. Meanwhile, CPTNN exhibits the optimal computational efficiency among all tested methods, with a computational time of only  $4.87$ seconds.

\begin{figure}[htbp]
    \centering
    \includegraphics[width=0.9\textwidth]{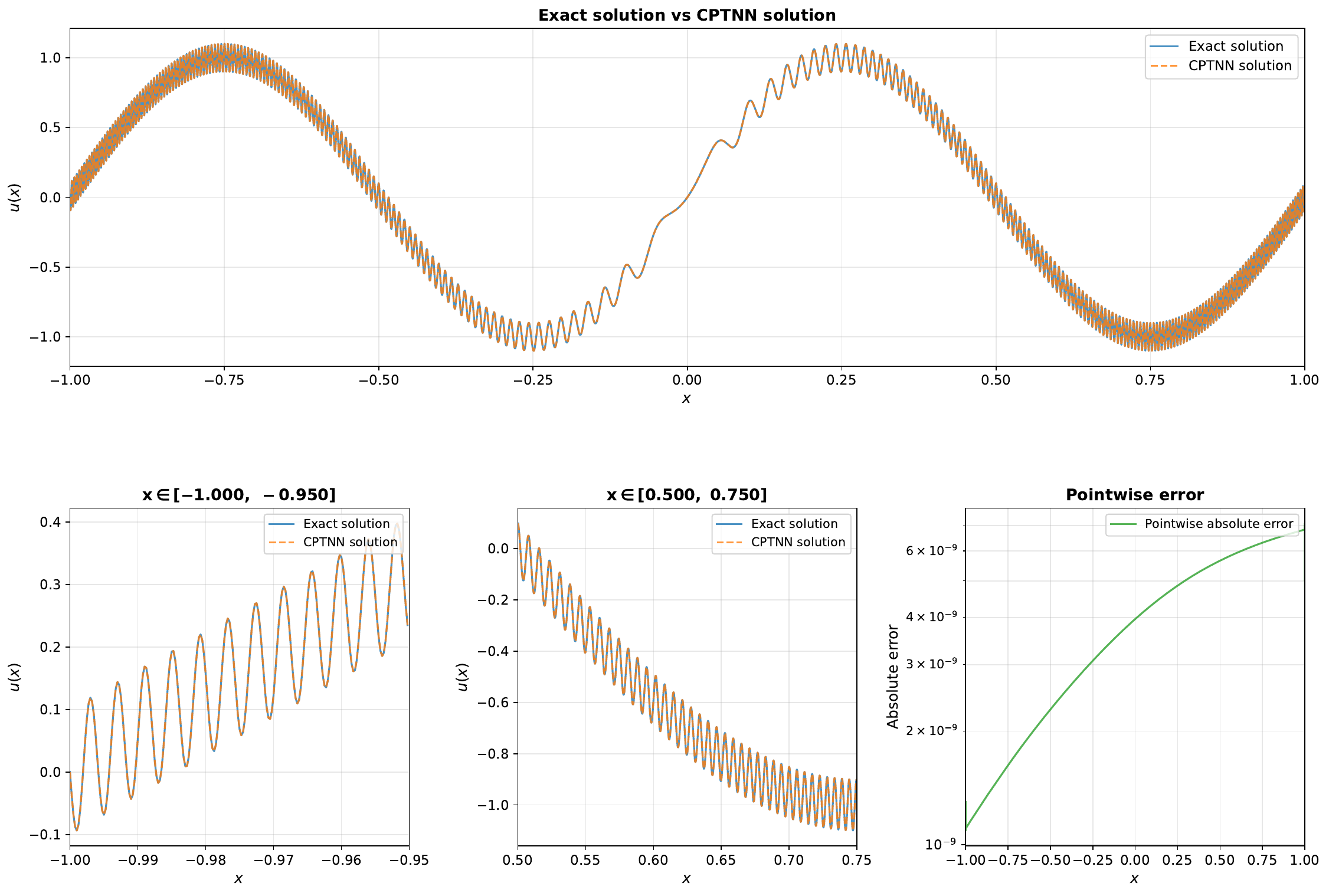}
    \caption{Numerical solution and absolute error obtained by CPTNN for equation \eqref{eq:variable coefficient}.}
    \label{fig:1D Poisson CPTNN}
\end{figure}

\begin{table}[!h]
    \centering
    \caption{Performance comparison of various methods for solving the equation \eqref{eq:variable coefficient} ($a=250$).}
    \label{tab:Variable coefficient Error&time}
    \begin{tabular}{c|c|c|c|c}
        \hline
        Method & PINN & RFM & TransNet & CPTNN \\
        \hline
        Error & 4.55E+00 & 4.96E-05 & 3.96E-06 & 6.13E-09 \\
        \hline
        Time & 878.60 & 22.89 & 22.29 & 4.87 \\
        \hline
    \end{tabular}
\end{table}

\subsubsection{Helmholtz equation}
Next, we consider the Helmholtz equation with constant wave numbers suggested in \cite{caiPhaseShiftDeep2020}:
\begin{equation}
\begin{aligned}
u''(x) + \lambda^2 u(x) = \mu^2 J_0''(\mu x) + \lambda^2 J_0(\mu x), x\in[-1,1],\quad\quad\\ 
u(-1) = J_0(-\mu)+0.2\cos(-\lambda),
u(1) = J_0(\mu)+0.2\cos(\lambda).
\end{aligned}
\label{eq:Helmholtz equation}
\end{equation}
where \(J_0(x)\) denotes  the zeroth-order Bessel function. The exact solution  $u(x)$ to this problem is given by:
\begin{equation}
u(x) = J_0(\mu x) + 0.2cos(\lambda x), 
\label{eq:Helmholtz solution}
\end{equation}

\begin{figure}[htbp]
    \centering
    \includegraphics[width=0.9\textwidth]{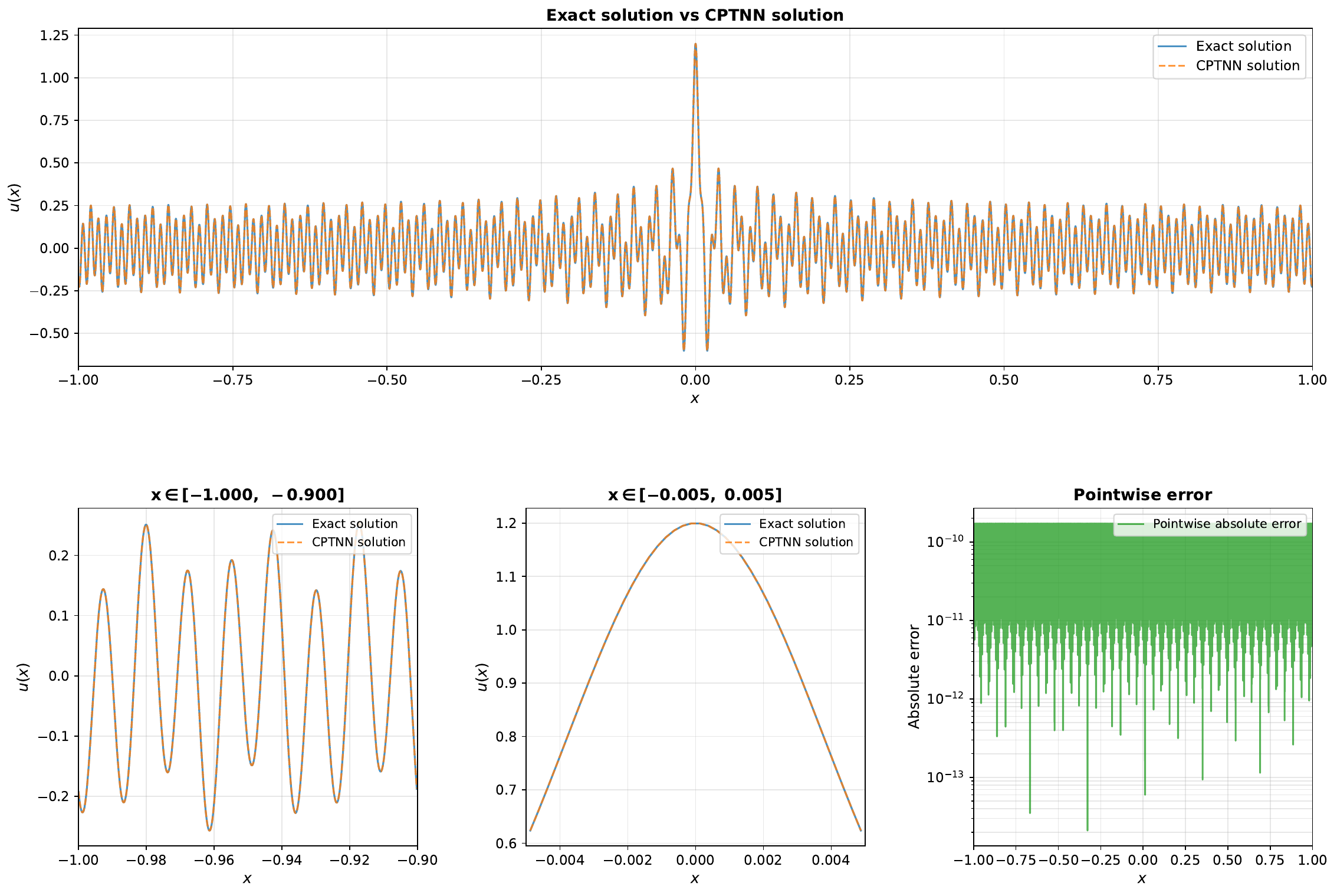}
    \caption{Numerical solution and absolute error obtained by CPTNN for Helmholtz equation \eqref{eq:Helmholtz equation}}
    \label{fig:1D Helmholtz CPTNN}
\end{figure}

In this numerical experiment, we set \(\lambda=500\) and \(\mu=200\). The number of hidden neurons is set to $1000$ for both PINN and CPTNN. For CPTNN, $500$ distinct frequencies are uniformly sampled from the interval $[0,200]$. For the baseline methods, we perform the following parameter settings:

\begin{itemize}
    \item For RFM, the scaling parameter $R_{\max}$ is optimized over the set 
    $\{1,11,\dots,351\}$.
    \item For TransNet, the shape parameter is exhaustively searched over the set $\{2,4,\dots,350\}$.
\end{itemize}
Figure \ref{fig:1D Helmholtz CPTNN} depicts the numerical solution and absolute error obtained by CPTNN for the large-wavenumber Helmholtz equation \eqref{eq:Helmholtz equation}, which clearly demonstrates the superior performance of the proposed method on this challenging high-frequency problem.

Table \ref{tab:Helmholtz Error&time} summarizes the relative $L_2$ errors and computational times of various methods for problem \eqref{eq:Helmholtz equation}.  
As evident from the table, the standard PINN fails completely due to spectral bias. Although RFM and TransNet are capable of capturing the oscillatory behavior of the solution, their accuracy saturates at the $\mathcal{O}(10^{-6})-\mathcal{O}(10^{-7})$ level. In contrast, the proposed CPTNN overcomes this accuracy bottleneck, reducing the relative $L_{2}$ error to $\mathcal{O}(10^{-10})$. Concurrently, our method exhibits superior computational efficiency, requiring significantly less time than all competing approaches considered.

\begin{table}[!h]
    \centering
    \caption{Performance comparison of various methods for solving the Helmholtz equation.}
    \label{tab:Helmholtz Error&time}
    \begin{tabular}{c|c|c|c|c}
        \hline
        Method & PINN & RFM & TransNet & CPTNN \\
        \hline
        Error & 9.99E-01 & 3.73E-06 & 2.15E-07 & 6.72E-10 \\
        \hline
        Time & 212.93 & 20.01 & 23.04 & 1.11 \\
        \hline
    \end{tabular}
\end{table}

\subsubsection{Nonlinear Helmholtz equation}
We next evaluate the CPTNN method on a nonlinear Helmholtz boundary value problem \cite{dongLocalExtremeLearning2021}: 
\begin{equation}
\begin{aligned}
u'' - 50 u + 10 \sin (u) = f(x) \\
u(-1) = h_{1}, \quad u(1) = h_{2},
\end{aligned}
\label{eq:nonlinear Helmholtz equation}
\end{equation}
with exact solution
\begin{equation}
u(x) = \sin\left(150\pi x + \frac{3\pi}{20}\right) \cos\left(200\pi x - \frac{2\pi}{5}\right) + \frac{3}{2} + \frac{x}{10},
\label{eq:nonlinear u(x)}
\end{equation}
where the source term $f(x)$ and the boundary values $h_1$ and $h_2$ are chosen to satisfy Eq. (\ref{eq:nonlinear Helmholtz equation}).

\begin{figure}[htbp]
    \centering
    \includegraphics[width=0.9\textwidth]{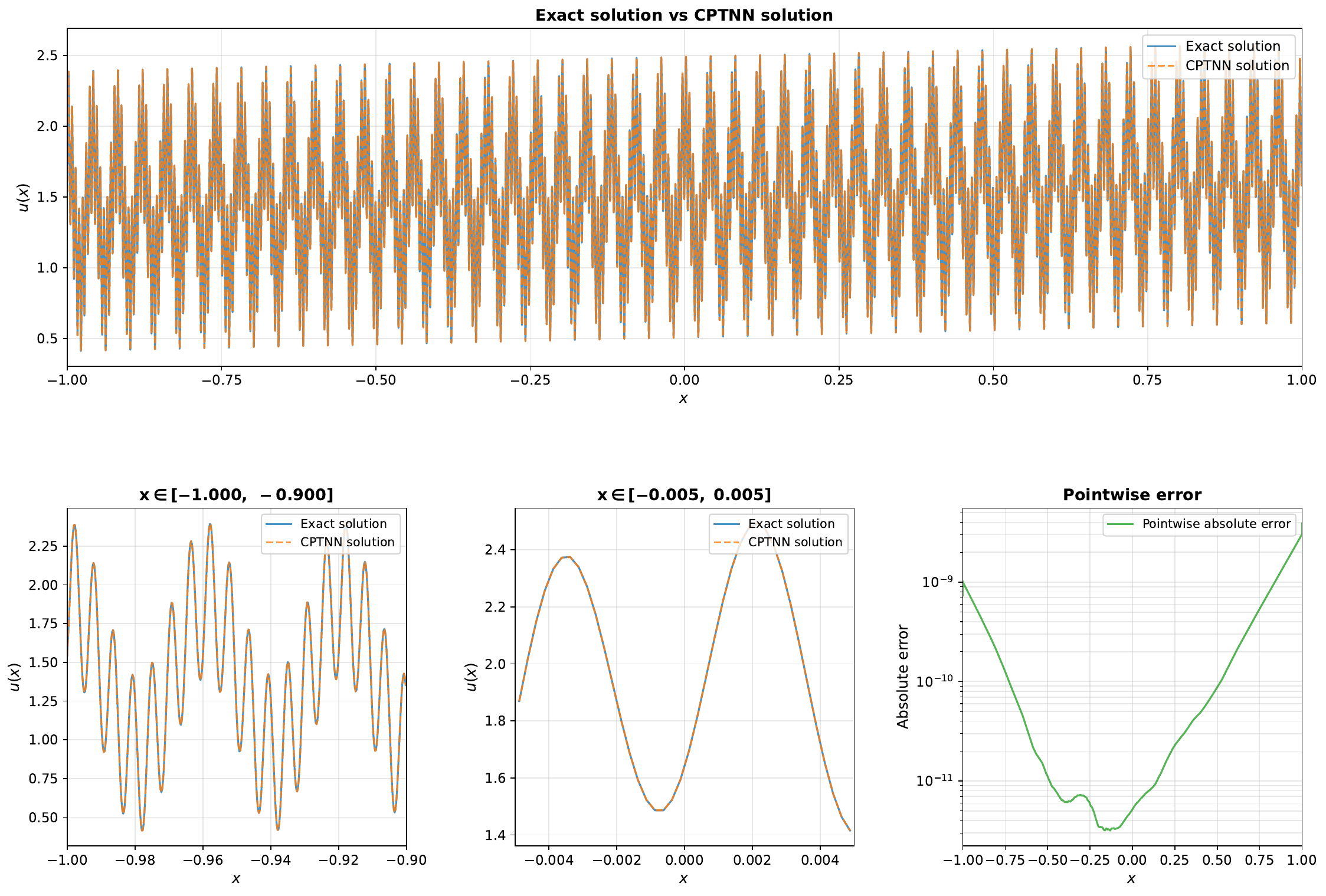}
    \caption{Numerical solution and absolute error obtained by CPTNN for nonlinear Helmholtz equation \eqref{eq:nonlinear Helmholtz equation}}
    \label{fig:1D Nonlinear Helmholtz CPTNN}
\end{figure}

For CPTNN,  we set $2000$ hidden neurons and sample $500$ distinct frequencies uniformly from interval $[0,250]$. For the baseline methods, we perform the following parameter tuning:
\begin{itemize}
    \item For RFM, the scaling parameter $R_{\max}$ is optimized over the set $\{1,11,\dots,201\}$.
    \item For TransNet, the shape parameter is exhaustively searched over the set $\{2,4,\dots,350\}$.
\end{itemize}
Figure \ref{fig:1D Nonlinear Helmholtz CPTNN} plots the numerical solution and absolute error obtained by CPTNN for the nonlinear large-wavenumber Helmholtz equation \eqref{eq:nonlinear Helmholtz equation}, demonstrating its excellent performance in capturing highly oscillatory nonlinear behavior.

Table \ref{tab:Nonlinear Helmholtz Error&time} summarizes the performance of various methods for solving the problem \eqref{eq:nonlinear Helmholtz equation}. The coupling of nonlinear terms and high-frequency oscillations significantly increases the problem difficulty: PINN fails to converge. While RFM and TransNet mitigate this challenge but incur high computational cost ($>250$s) due to over-parameterization and exhaustive tuning, with accuracy saturating at
$\mathcal{O}(10^{-7})$. In stark contrast, CPTNN achieves a relative $L_2$ error of $\mathcal{O}(10^{-10})$ (\textbf{$3$ orders of magnitude} higher accuracy than the best-performing baseline) with only a fraction of the computational cost, showcasing its exceptional accuracy and efficiency.

\begin{table}[!h]
    \centering
    \caption{Performance comparison of various methods for solving the problem \eqref{eq:nonlinear Helmholtz equation}}
    \label{tab:Nonlinear Helmholtz Error&time}
    \begin{tabular}{c|c|c|c|c}
        \hline
        Method & PINN & RFM & TransNet & CPTNN \\
        \hline
        Error & 2.04E+01 & 5.19E-07 & 1.58E-07 & 3.84E-10 \\
        \hline
        Time & 380.99 & 250.70 & 258.98 & 35.87 \\
        \hline
    \end{tabular}
\end{table}

\subsubsection{Wave equation}
\begin{figure}[htbp]
    \centering
    \includegraphics[width=0.9\textwidth]{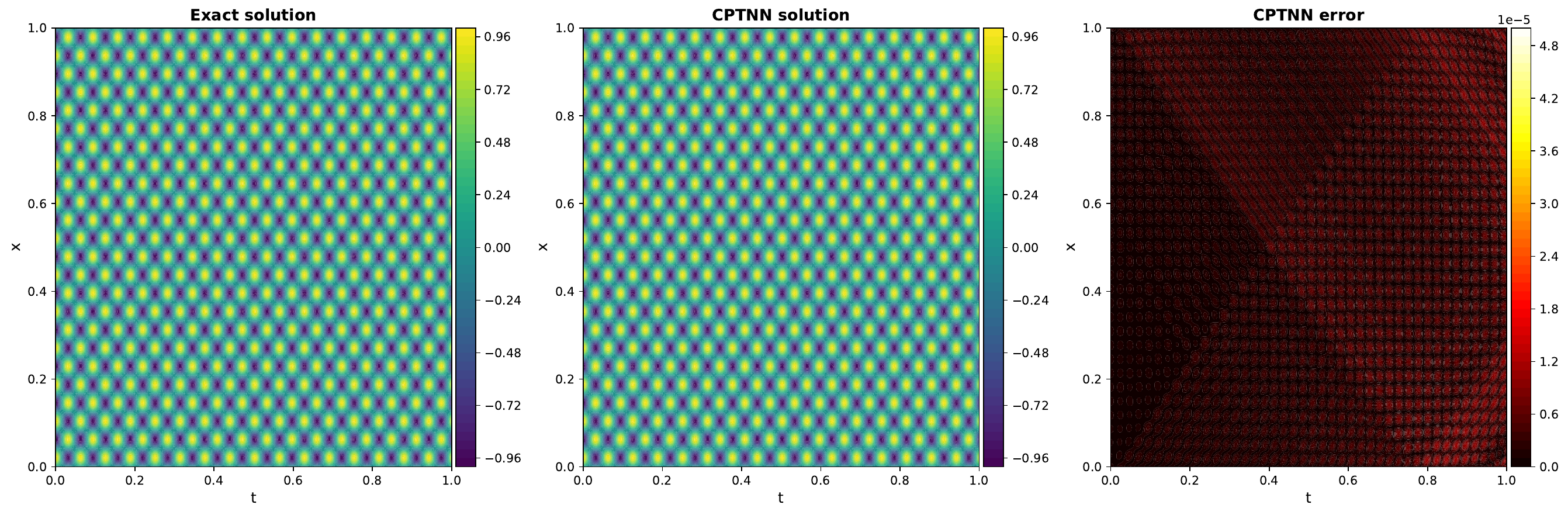}
    \caption{Numerical solution and absolute error obtained by CPTNN for Wave equation \eqref{eq:Wave equation}.}
    \label{fig:2d wave}
\end{figure}

Here, we consider the following wave equation \cite{zhangTransferableNeuralNetworks2024}:
\begin{equation}
\begin{aligned}
&u_{tt}(x,t)-cu_{xx}(x,t)=0, && (x,t)\in[0,1]\times[0,1], \\
&u(x,0)=\sin(24\pi x), && x\in[0,1], \\
&u_{t}(x,0)=0, && x\in[0,1], \\
&u(0,t)=u(1,t)=0, && t\in[0,1],
\end{aligned}
\label{eq:Wave equation}
\end{equation}
where $c=\frac{625}{36 \pi^2}$. The exact solution is given by
\begin{equation}
u(x,t) = 0.5 \left( \sin(24 \pi x + 100 t) + \sin(24 \pi x - 100 t) \right).
\label{eq:2d wave}
\end{equation}
For PINN and CPTNN, we fix the number of hidden neurons at $4000$. For CPTNN, we select  $20^2$ distinct frequencies uniformly sampled from $[-20, 20]^2$. For the baseline methods, we perform the following parameter tuning:
\begin{itemize}
    \item For RFM, the scaling parameter $R_{\max}$ is optimized over the set $\{1,11,\dots,101\}$.
    \item For TransNet, the shape parameter is exhaustively searched over the set $\{2,4,\dots,30\}$.
\end{itemize}
Figure \ref{fig:2d wave} depicts the numerical solution and absolute error obtained by CPTNN for the wave equation \eqref{eq:Wave equation}, demonstrating its excellent performance in capturing highly oscillatory solutions for time-dependent wave problems.

Table \ref{tab:Wave Error&time} summarizes the relative $L_2$ errors and computational costs for solving this time-dependent wave equation. The standard PINN fails to produce a valid solution. While RFM and TransNet successfully solve the equation, they incur significant computational overhead (approximately 500 seconds) due to the excessive number of hidden neurons. In contrast, the proposed CPTNN achieves a relative $L_2$ errors of  $\mathcal{O}(10^{-6})$, outperforming the best baseline method by nearly two orders of magnitude. Furthermore, CPTNN obtains this high-fidelity solution in only 15.81 seconds, demonstrating its superior computational efficiency.

\begin{table}[!h]
    \centering
    \caption{Performance comparison of various methods for solving the Wave equation.}
    \label{tab:Wave Error&time}
    \begin{tabular}{c|c|c|c|c}
        \hline
        Method & PINN & RFM & TransNet & CPTNN \\
        \hline
        Error & 1.00E+00 & 5.40E-03 & 2.23E-04 & 5.95E-06 \\
        \hline
        Time & 13213.59 & 580.63 & 496.22 & 15.81 \\
        \hline
    \end{tabular}
\end{table}

\subsubsection{2D Diffusion interface problem}
In the final numerical example, we consider the following diffusion interface problem  \cite{luMultipleTransferableNeural2025} posed on the domain $\Omega = [0,2]^2$:
\begin{equation}
\begin{aligned}
-\nabla \cdot (\beta \nabla u) &= f, \quad (x, y) \in \Omega_1 \cup \Omega_2, \\
[u] &= h_1, \quad (x, y) \in \Gamma, \\
[\beta \nabla u \cdot n] &= h_2, \quad (x, y) \in \Gamma, \\
u &= g, \quad (x, y) \in \partial\Omega,
\label{eq:Interface-problem}
\end{aligned}
\end{equation}
where $\Gamma = \{(x, y) \mid (x - 1)^2 + (y - 1)^2 = 0.5\}$  denotes the circular interface separating the subdomains $\Omega_1$ and  $\Omega_2$, and the diffusion coefficient is given by
$\beta(x, y) = 1 \cdot \chi_{\Omega_1} + 10 \cdot \chi_{\Omega_2}$. The exact solution is given by

\begin{equation}
u(x,y) =
\begin{cases}
(x-1)^2 + (y-1)^2, & (x,y) \in \Omega_1, \\
\cos(16\pi x) \cos(16\pi y), & (x,y) \in \Omega_2.
\end{cases}
\end{equation}
The solution to problem \eqref{eq:Interface-problem} exhibits distinct frequency characteristics across the two subdomains: $\Omega_1$ admits a smooth, low-frequency solution, while $\Omega_2$ features a highly oscillatory high-frequency component. 

To address the multiscale discrepancies inherent to the interface problem, we employ the domain decomposition strategy from \cite{luMultipleTransferableNeural2025} and incorporate a dual-subnetwork framework. For consistency and fair comparison, this strategy is applied to all baseline models.
The number of hidden neurons per subnetwork is set to 5000 for PINN and CPTNN, and 50000 for RFM and TransNet. For the baseline methods, we tune parameters as follows:
\begin{itemize}
    \item For RFM, $R_{\max,1} \in \{1,11,\dots,51\}$, $R_{\max,2} \in \{1,11,\dots,101\}$.
    \item For TransNet, $\gamma_1 \in \{2,4,\dots,10\}, \gamma_2 \in \{2,4,\dots,30\}$.
\end{itemize}
In stark contrast, both subnetworks in the CPTNN framework share \textbf{identical hyperparameters} throughout training, removing the burden of manual tuning entirely. For frequency sampling, $50^2$ distinct modes are uniformly sampled over $[-10, 10]^2$. Interior training/testing points follow the setup from prior examples, with $100$ collocation points on each edge of $\partial\Omega$ and $300$ on $\Gamma$ to enforce interface jump conditions.

\begin{figure}[htbp]
    \centering
    \includegraphics[width=0.9\textwidth]{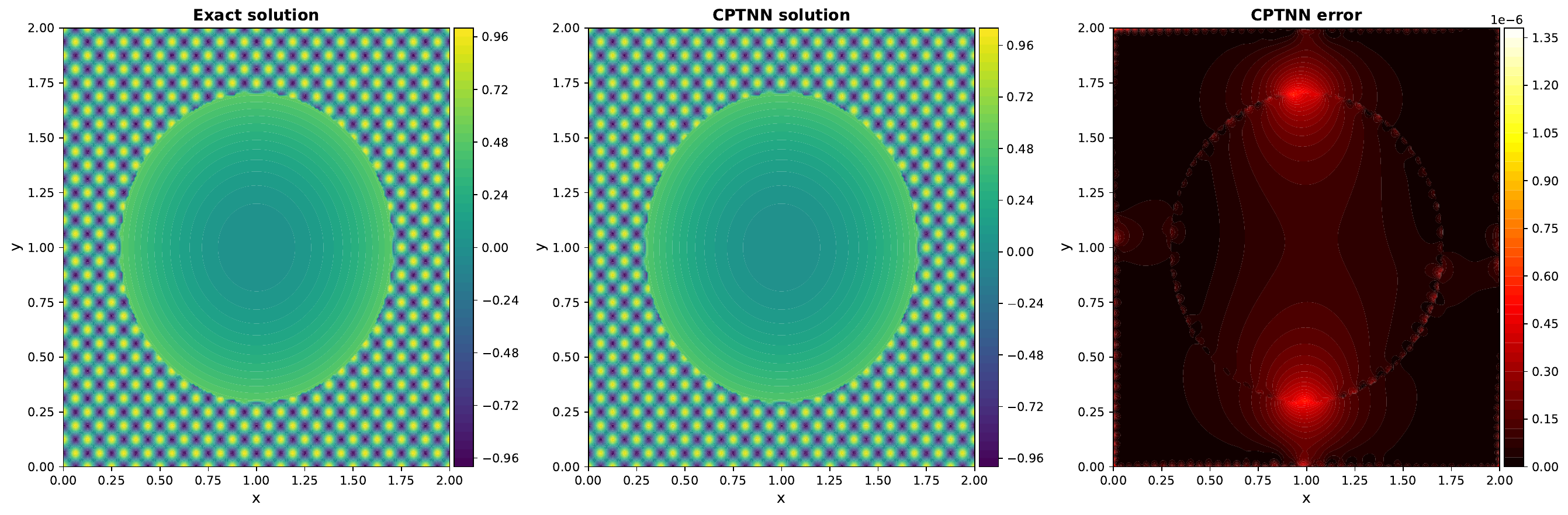}
    \caption{Numerical solution and absolute error obtained by CPTNN for problem \eqref{eq:Interface-problem}.}
    \label{fig:2d interface}
\end{figure}

Figure \ref{fig:2d interface} shows the numerical solution and absolute error of CPTNN for the problem \eqref{eq:Interface-problem}, confirming the method’s successful extension to interface problems with a maximum absolute error of $\mathcal{O}(10^{-6})$.

Table \ref{tab:Interface Error&time} compares the relative $L_2$ errors and computational times of numerical solutions obtained by various methods for solving Eq. \eqref{eq:Interface-problem}. Despite utilizing identical hyperparameters across both sub-networks and avoiding manual tuning, the proposed CPTNN achieves a relative $L_2$ error of $\mathcal{O}(10^{-7})$, outperforming the best baseline method by three orders of magnitude with minimal computational cost.

\begin{table}[!h]
    \centering
    \caption{Performance comparison of various methods for problem \eqref{eq:Interface-problem}.}
    \label{tab:Interface Error&time}
    \begin{tabular}{c|c|c|c|c}
        \hline
        Method & PINN & RFM & TransNet & CPTNN \\
        \hline
        Error & 1.71E+02 & 1.42E-03 & 9.11E-04 & 2.32E-07 \\
        \hline
        Time & 13819.51 & 507.86 & 456.15 & 113.54 \\
        \hline
    \end{tabular}
\end{table}

\section{Concluding remarks}
\label{sec:conclusion}
In this work, we propose the \textbf{PhaseTNN} framework to address critical bottlenecks in neural network-based high-frequency function approximation and PDE solving: spectral bias and the accuracy-efficiency trade-off induced by empirical parameter tuning. PhaseTNN integrates improved frequency filtering, phase-shift frequency decomposition, and TransNet’s efficient least-squares training to decompose intractable high-frequency problems into tractable low-frequency subproblems, eliminating the need for manual shape parameter tuning. We present two PhaseTNN architectures:
\begin{itemize}
    \item \textbf{PPTNN}: Decomposes high-frequency functions into narrow-band components for parallel TransNet training, achieving high-precision 1D/2D function approximation;
    \item \textbf{CPTNN}: Integrates phase-modulated basis functions into a single model, enabling simultaneous multi-frequency representation and efficient solution of linear/nonlinear high-frequency PDEs, with scalability to high-dimensional problems.
\end{itemize}
Extensive numerical experiments validate that PhaseTNN achieves near-machine-precision accuracy for 1D/2D high-frequency function approximation, and outperforms state-of-the-art methods (PINN, RFM, TransNet, MSNN) in both accuracy and computational efficiency for canonical high-frequency PDEs (variable-coefficient elliptic equation, linear/nonlinear Helmholtz, wave equations, 2D diffusion interface problem).

Our future research directions include:
\begin{itemize}
 \item \textbf{Uncertainty quantification}: Integrate PhaseTNN with stochastic methods to solve stochastic PDEs with high-frequency random solutions;
 \item \textbf{Non-smooth functions}: Extend the filtering method to handle non-smooth high-frequency functions (e.g., discontinuous solutions) and singular PDEs;
 \item \textbf{Inverse problems}: Apply PhaseTNN to high-frequency PDE inverse problems (e.g., parameter identification) for scientific computing and engineering applications.
 \end{itemize}
\section*{Acknowledgments}
The work is supported by National Natural Science Foundation of China (No. 92370113, 12326346, 12461019),
Natural Science Foundation of Zhejiang Province (No. ZCLY24A0101), the ‘‘Pioneer’’ and ‘‘Leading Goose’’ R\&D Program of Zhejiang Province (No. 2024C1093), and the Science and Technology Research Project of Jiangxi Provincial Department of Education 
(No. GJJ2504901). 
%

\begin{thebibliography}{999}
\addtolength{\itemsep}{-0.5ex}

\bibitem{ainsworthGalerkinNeuralNetwork2022}
Ainsworth, M., Dong, J. Galerkin neural network approximation of singularly-perturbed elliptic systems. \textit{Computer Methods in Applied Mechanics and Engineering}, 2022, 402: 115169.

\bibitem{baldiNeuralNetworksPrincipal1989}
Baldi, P., Hornik, K. Neural networks and principal component analysis: learning from examples without local minima. \textit{Neural Networks}, 1989, 2(1): 53--58.

\bibitem{bottouOptimizationMethodsLargescale2018}
Bottou, L., Curtis, F.E., Nocedal, J. Optimization methods for large-scale machine learning. \textit{SIAM Review}, 2018, 60(2): 223--311.

\bibitem{brownLanguageModelsAre2020}
Brown, T., et al. Language models are few-shot learners. \textit{Advances in Neural Information Processing Systems}, 2020, 33: 1877--1901.

\bibitem{brownStatisticalApproachMachine1990a}
Brown, P.F., et al. A statistical approach to machine translation. \textit{Computational Linguistics}, 1990, 16(2): 79--85.

\bibitem{caiPhaseShiftDeep2020}
Cai, W., Li, X., Liu, L. A phase shift deep neural network for high frequency approximation and wave problems. \textit{SIAM Journal on Scientific Computing}, 2020, 42(5): A3285--A3312.

\bibitem{chenBridgingTraditionalMachine2022}
Chen, J., Chi, X., E, W., Yang, Z. Bridging traditional and machine learning-based algorithms for solving PDEs: the random feature method. \textit{Journal of Machine Learning}, 2022, 1(3): 268--298.

\bibitem{chenRandomFeatureMethod2023}
Chen, J.-R., E, W., Luo, Y.-X. The random feature method for time-dependent problems. \textit{East Asian Journal on Applied Mathematics}, 2023, 13(3): 435--463.

\bibitem{chiuCANPINNFastPhysicsinformed2022}
Chiu, P.-H., et al. CAN-PINN: a fast physics-informed neural network based on coupled automatic--numerical differentiation method. \textit{Computer Methods in Applied Mechanics and Engineering}, 2022, 395: 114909.

\bibitem{chowdharyNaturalLanguageProcessing2020}
Chowdhary, K.R. Natural language processing. In: \textit{Fundamentals of Artificial Intelligence}, Springer, 2020: 603--649.

\bibitem{collobertUnifiedArchitectureNatural2008}
Collobert, R., Weston, J. A unified architecture for natural language processing: deep neural networks with multitask learning. In: \textit{ICML}, 2008: 160--167.

\bibitem{devlinBertPretrainingDeep2019}
Devlin, J., Chang, M.-W., Lee, K., Toutanova, K. BERT: pre-training of deep bidirectional transformers for language understanding. In: \textit{NAACL-HLT}, 2019: 4171--4186.

\bibitem{dongComputingHyperparameterExtreme2022}
Dong, S., Yang, J. On computing the hyperparameter of extreme learning machines: algorithm and application to computational PDEs, and comparison with classical and high-order finite elements. \textit{Journal of Computational Physics}, 2022, 463: 111290.

\bibitem{dongLocalExtremeLearning2021}
Dong, S., Li, Z. Local extreme learning machines and domain decomposition for solving linear and nonlinear partial differential equations. \textit{Computer Methods in Applied Mechanics and Engineering}, 2021, 387: 114129.

\bibitem{dwivediPhysicsInformedExtreme2020}
Dwivedi, V., Srinivasan, B. Physics informed extreme learning machine (PIELM)--a rapid method for the numerical solution of partial differential equations. \textit{Neurocomputing}, 2020, 391: 96--118.

\bibitem{hanDeepLearningbasedNumerical2017}
Han, J., Jentzen, A. Deep learning-based numerical methods for high-dimensional parabolic partial differential equations and backward stochastic differential equations. \textit{Communications in Mathematics and Statistics}, 2017, 5(4): 349--380.

\bibitem{hanSolvingHighdimensionalPartial2018}
Han, J., Jentzen, A., E, W. Solving high-dimensional partial differential equations using deep learning. \textit{PNAS}, 2018, 115(34): 8505--8510.

\bibitem{heDeepResidualLearning2016}
He, K., Zhang, X., Ren, S., Sun, J. Deep residual learning for image recognition. In: \textit{CVPR}, 2016: 770--778.

\bibitem{hornikMultilayerFeedforwardNetworks1989}
Hornik, K., Stinchcombe, M., White, H. Multilayer feedforward networks are universal approximators. \textit{Neural Networks}, 1989, 2(5): 359--366.

\bibitem{howardStackedNetworksImprove2023}
Howard, A.A., et al. Stacked networks improve physics-informed training: applications to neural networks and deep operator networks. arXiv preprint arXiv:2311.06483, 2023.

\bibitem{huangExtremeLearningMachine2006}
Huang, G.-B., Zhu, Q.-Y., Siew, C.-K. Extreme learning machine: theory and applications. \textit{Neurocomputing}, 2006, 70(1--3): 489--501.

\bibitem{jagtapExtendedPhysicsinformedNeural2020}
Jagtap, A.D., Karniadakis, G.E. Extended physics-informed neural networks (XPINNs): a generalized space-time domain decomposition based deep learning framework for nonlinear partial differential equations. \textit{Communications in Computational Physics}, 2020, 28(5).

\bibitem{johnxuFrequencyPrincipleFourier2020}
John Xu, Z.-Q., Zhang, Y., Luo, T., Xiao, Y., Ma, Z. Frequency principle: Fourier analysis sheds light on deep neural networks. \textit{Communications in Computational Physics}, 2020, 28(5): 1746--1767.

\bibitem{kiyaniOptimizingOptimizerPhysicsinformed2025}
Kiyani, E., et al. Optimizing the optimizer for physics-informed neural networks and Kolmogorov-Arnold networks. \textit{Computer Methods in Applied Mechanics and Engineering}, 2025, 446: 118308.

\bibitem{kopanicakovaEnhancingTrainingPhysicsinformed2024}
Kopaničáková, A., Kothari, H., Karniadakis, G.E., Krause, R. Enhancing training of physics-informed neural networks using domain decomposition--based preconditioning strategies. \textit{SIAM Journal on Scientific Computing}, 2024, 46(5): S46--S67.

\bibitem{krishnapriyanCharacterizingPossibleFailure2021}
Krishnapriyan, A., et al. Characterizing possible failure modes in physics-informed neural networks. \textit{Advances in Neural Information Processing Systems}, 2021, 34: 26548--26560.

\bibitem{lecunDeepLearning2015}
LeCun, Y., Bengio, Y., Hinton, G. Deep learning. \textit{Nature}, 2015, 521(7553): 436--444.

\bibitem{liuKanKolmogorovarnoldNetworks2024}
Liu, Z., et al. KAN: Kolmogorov-Arnold networks. arXiv preprint arXiv:2404.19756, 2024.

\bibitem{luMultipleTransferableNeural2025}
Lu, T., Ju, L., Zhu, L. A multiple transferable neural network method with domain decomposition for elliptic interface problems. \textit{Journal of Computational Physics}, 2025, 530: 113902.

\bibitem{luoUpperLimitDecaying2022}
Luo, T., Ma, Z., Wang, Z., Xu, Z.J., Zhang, Y. An upper limit of decaying rate with respect to frequency in linear frequency principle model. In: \textit{MSML}, PMLR, 2022: 205--214.

\bibitem{mildenhallNerfRepresentingScenes2021}
Mildenhall, B., et al. NeRF: representing scenes as neural radiance fields for view synthesis. \textit{Communications of the ACM}, 2021, 65(1): 99--106.

\bibitem{moseleyFiniteBasisPhysicsinformed2023}
Moseley, B., Markham, A., Nissen-Meyer, T. Finite basis physics-informed neural networks (FBPINNs): a scalable domain decomposition approach for solving differential equations. \textit{Advances in Computational Mathematics}, 2023, 49(4): 62.

\bibitem{mullerAchievingHighAccuracy2023}
Müller, J., Zeinhofer, M. Achieving high accuracy with PINNs via energy natural gradient descent. In: \textit{ICML}, PMLR, 2023: 25471--25485.

\bibitem{niNumericalComputationPartial2023}
Ni, N., Dong, S. Numerical computation of partial differential equations by hidden-layer concatenated extreme learning machine. \textit{Journal of Scientific Computing}, 2023, 95(2): 35.

\bibitem{rahamanSpectralBiasNeural2019}
Rahaman, N., et al. On the spectral bias of neural networks. In: \textit{ICML}, PMLR, 2019: 5301--5310.

\bibitem{raissiPhysicsinformedNeuralNetworks2019}
Raissi, M., Perdikaris, P., Karniadakis, G.E. Physics-informed neural networks: a deep learning framework for solving forward and inverse problems involving nonlinear partial differential equations. \textit{Journal of Computational Physics}, 2019, 378: 686--707.

\bibitem{ronnebergerUnetConvolutionalNetworks2015}
Ronneberger, O., Fischer, P., Brox, T. U-Net: convolutional networks for biomedical image segmentation. In: \textit{MICCAI}, 2015: 234--241.

\bibitem{shuklaComprehensiveFAIRComparison2024}
Shukla, K., et al. A comprehensive and FAIR comparison between MLP and KAN representations for differential equations and operator networks. \textit{Computer Methods in Applied Mechanics and Engineering}, 2024, 431: 117290.

\bibitem{sirignanoDGMDeepLearning2018}
Sirignano, J., Spiliopoulos, K. DGM: a deep learning algorithm for solving partial differential equations. \textit{Journal of Computational Physics}, 2018, 375: 1339--1364.

\bibitem{urbanUnveilingOptimizationProcess2025}
Urbán, J.F., Stefanou, P., Pons, J.A. Unveiling the optimization process of physics informed neural networks: how accurate and competitive can PINNs be? \textit{Journal of Computational Physics}, 2025, 523: 113656.

\bibitem{vaswaniAttentionAllYou2017}
Vaswani, A., et al. Attention is all you need. \textit{Advances in Neural Information Processing Systems}, 2017, 30.

\bibitem{wangGradientAlignmentPhysicsinformed2025}
Wang, S., Bhartari, A.K., Li, B., Perdikaris, P. Gradient alignment in physics-informed neural networks: a second-order optimization perspective. arXiv preprint arXiv:2502.00604, 2025.

\bibitem{wangMultistageNeuralNetworks2024}
Wang, Y., Lai, C.-Y. Multi-stage neural networks: function approximator of machine precision. \textit{Journal of Computational Physics}, 2024, 504: 112865.

\bibitem{wangUnderstandingMitigatingGradient2021}
Wang, S., Teng, Y., Perdikaris, P. Understanding and mitigating gradient flow pathologies in physics-informed neural networks. \textit{SIAM Journal on Scientific Computing}, 2021, 43(5): A3055--A3081.

\bibitem{wangWhenWhyPINNs2022}
Wang, S., Yu, X., Perdikaris, P. When and why PINNs fail to train: a neural tangent kernel perspective. \textit{Journal of Computational Physics}, 2022, 449: 110768.

\bibitem{xuMultigradeDeepLearning2025}
Xu, Y. Multi-grade deep learning. \textit{Communications on Applied Mathematics and Computation}, 2025: 1--52.

\bibitem{xuNovelClassHessian2025}
Xu, M., Zhang, L., Wu, B., Liu, K. A novel class of hessian recovery-based numerical methods for solving biharmonic equations and their applications in phase field modeling. \textit{Finite Elements in Analysis and Design}, 2025, 251: 104405.

\bibitem{xuOverviewFrequencyPrinciple2025}
Xu, Z.-Q.J., Zhang, Y., Luo, T. Overview frequency principle/spectral bias in deep learning. \textit{Communications on Applied Mathematics and Computation}, 2025, 7(3): 827--864.

\bibitem{xuUnderstandingOvercomingSpectral2025a}
Xu, Z.-Q.J., Zhang, L., Cai, W. On understanding and overcoming spectral biases of deep neural network learning methods for solving PDEs. \textit{Journal of Computational Physics}, 2025, 530: 113905.

\bibitem{yangNovelImprovedExtreme2018}
Yang, Y., Hou, M., Luo, J. A novel improved extreme learning machine algorithm in solving ordinary differential equations by Legendre neural network methods. \textit{Advances in Difference Equations}, 2018, 2018(1): 469.

\bibitem{yuDeepRitzMethod2018}
Yu, B. The deep Ritz method: a deep learning-based numerical algorithm for solving variational problems. \textit{Communications in Mathematics and Statistics}, 2018, 6(1): 1--12.

\bibitem{zangWeakAdversarialNetworks2020}
Zang, Y., Bao, G., Ye, X., Zhou, H. Weak adversarial networks for high-dimensional partial differential equations. \textit{Journal of Computational Physics}, 2020, 411: 109409.

\bibitem{zhangNewFiniteElement2005}
Zhang, Z., Naga, A. A new finite element gradient recovery method: superconvergence property. \textit{SIAM Journal on Scientific Computing}, 2005, 26(4): 1192--1213.

\bibitem{zhangTransferableNeuralNetworks2024}
Zhang, Z., Bao, F., Ju, L., Zhang, G. Transferable neural networks for partial differential equations. \textit{Journal of Scientific Computing}, 2024, 99(1): 2.

\end{thebibliography}

\end{document}